\documentclass[11pt]{article}

\usepackage{amsmath, amsthm, amssymb}
\usepackage{enumitem}
\usepackage{pdflscape}
\usepackage{caption}
\usepackage{bm}
\usepackage{ifpdf}
	\ifpdf
\usepackage[pdftex]{graphicx}
	\else
\usepackage[dvips]{graphicx}
	\fi
\usepackage[all]{xy}
\usepackage{tocvsec2}
\usepackage{bbm}
	\input xy
	\xyoption{all}
\usepackage[pdftex,plainpages=false,hypertexnames=false,pdfpagelabels]{hyperref}
	\newcommand{\arxiv}[1]{\href{http://arxiv.org/abs/#1}{\tt arXiv:\nolinkurl{#1}}}
	\newcommand{\arXiv}[1]{\href{http://arxiv.org/abs/#1}{\tt arXiv:\nolinkurl{#1}}}
	\newcommand{\doi}[1]{\href{http://dx.doi.org/#1}{{\tt DOI:#1}}}

	\newcommand{\googlebooks}[1]{(preview at \href{http://books.google.com/books?id=#1}{google books})}
	
\usepackage{xcolor}
	\definecolor{dark-red}{rgb}{0.7,0.25,0.25}
	\definecolor{dark-blue}{rgb}{0.15,0.15,0.55}
	\definecolor{medium-blue}{rgb}{0,0,.8}
	\definecolor{DarkGreen}{RGB}{0,150,0}
	\definecolor{rho}{named}{red}
	\hypersetup{
	   colorlinks, linkcolor={purple},
	   citecolor={medium-blue}, urlcolor={medium-blue}
	}
\usepackage{longtable}
\usepackage{fullpage}

\usepackage{mathbbol}

\theoremstyle{plain}
\newtheorem{thm}{Theorem}[section]
\newtheorem*{thm*}{Theorem}
\newtheorem{thmalpha}{Theorem}

\newtheorem{cor}[thm]{Corollary}

\newtheorem*{cor*}{Corollary}

\newtheorem*{conj*}{Conjecture}
\newtheorem{lem}[thm]{Lemma}
\newtheorem{prop}[thm]{Proposition}

\newtheorem*{quest*}{Question}
\newtheorem*{claim*}{Claim}

\theoremstyle{definition}
\newtheorem{defn}[thm]{Definition}

\newtheorem{exs}[thm]{Examples}
\newtheorem{ex}[thm]{Example}
\newtheorem*{ex*}{Example}
\newtheorem{sub-ex}[thm]{Sub-Example}
\newtheorem{counter-ex}[thm]{Counter-Example}

\newtheorem*{rem*}{Remark}
\newtheorem{remark}[thm]{Remark}

\DeclareMathOperator{\Aut}{Aut}

\DeclareMathOperator{\End}{End}

\DeclareMathOperator{\op}{op}

\DeclareMathOperator{\supp}{supp}
\DeclareMathOperator{\id}{id}

\DeclareMathOperator{\Irr}{Irr}

\DeclareMathOperator{\Tr}{Tr}
\DeclareMathOperator{\tr}{tr}

\newcommand{\comment}[1]{}

\newcommand{\bbOne}{\mathbbm{1}}

\newcommand{\noshow}[1]{}
\newcommand{\MR}[1]{}

\newcommand{\Rep}{{\sf Rep}}

\newcommand{\fgpBim}{{\sf Bim_{fgp}}}
\newcommand{\fgpMod}{{\sf Mod_{fgp}}}

\newcommand{\spbfBim}{{\sf Bim_{bf}^{sp}}}
\newcommand{\spfgpBim}{{\sf Bim_{fgp}^{sp}}}
\renewcommand{\Vec}{{\sf Vec}}

\newcommand{\fdHilb}{{\sf Hilb_{fd}}}
\newcommand{\rCorr}{{\mathsf{C^{*}Alg}}}

\newcommand{\cCCAlgs}{{\sf C^*Alg}(\cC)}
\newcommand{\CDisc}{{\sf C^*Disc}}

\newcommand{\PQN}{{\sf PQN}}
\newcommand{\QN}{{\sf QN}}
\newcommand{\<}{\langle}
\renewcommand{\>}{\rangle}

\def\semicolon{;}
\def\applytolist#1{
    \expandafter\def\csname multi#1\endcsname##1{
        \def\multiack{##1}\ifx\multiack\semicolon
            \def\next{\relax}
        \else
            \csname #1\endcsname{##1}
            \def\next{\csname multi#1\endcsname}
        \fi
        \next}
    \csname multi#1\endcsname}

\def\calc#1{\expandafter\def\csname c#1\endcsname{{\mathcal #1}}}
\applytolist{calc}QWERTYUIOPLKJHGFDSAZXCVBNM;
\def\bbc#1{\expandafter\def\csname bb#1\endcsname{{\mathbb #1}}}
\applytolist{bbc}QWERTYUIOPLKJHGFDSAZXCVBNM;
\def\bfc#1{\expandafter\def\csname bf#1\endcsname{{\mathbf #1}}}
\applytolist{bfc}QWERTYUIOPLKJHGFDSAZXCVBNM;
\def\sfc#1{\expandafter\def\csname s#1\endcsname{{\sf #1}}}
\applytolist{sfc}QWERTYUIOPLKJHGFDSAZXCVBNM;
\def\fc#1{\expandafter\def\csname f#1\endcsname{{\mathfrak #1}}}
\applytolist{fc}QWERTYUIOPLKJHGFDSAZXCVBNM;

\usepackage{tikz}
\usepackage{tikz-cd}
\usetikzlibrary{arrows,backgrounds,patterns.meta}
\usetikzlibrary{positioning,shadings,cd}
\usetikzlibrary{shapes}
\usetikzlibrary{backgrounds}
\usetikzlibrary{decorations,decorations.pathreplacing,decorations.markings}
\usetikzlibrary{fit,calc,through}
\usetikzlibrary{external}
\usetikzlibrary{arrows}
\tikzset{vertex/.style = {shape=circle,draw,fill=black,inner sep=0pt,minimum size=5pt}}
\tikzset{edge/.style = {->,> = latex', bend right}}
\tikzset{
	super thick/.style={line width=3pt}
}
\tikzset{
    quadruple/.style args={[#1] in [#2] in [#3] in [#4]}{
        #1,preaction={preaction={preaction={draw,#4},draw,#3}, draw,#2}
    }
}
\tikzstyle{shaded}=[fill=red!10!blue!20!gray!30!white]
\tikzstyle{unshaded}=[fill=white]
\tikzstyle{empty box}=[circle, draw, thick, fill=white, opaque, inner sep=2mm]
\tikzstyle{annular}=[scale=.7, inner sep=1mm, baseline]
\tikzstyle{rectangular}=[scale=.75, inner sep=1mm, baseline=-.1cm]
\tikzstyle{mid>}=[decoration={markings, mark=at position 0.5 with {\arrow{>}}}, postaction={decorate}]
\tikzstyle{mid<}=[decoration={markings, mark=at position 0.5 with {\arrow{<}}}, postaction={decorate}]
\tikzstyle{over}=[double, draw=white, super thick, double=]

\tikzdeclarepattern{
  name=primeddots,
  type=uncolored,
  bounding box={(-.6pt,-.6pt) and (.6pt,.6pt)},
  tile size={(5pt,5pt)},
  tile transformation={rotate=60},
  parameters={none}, 
  code={
    \fill(0pt,0pt) circle (.35pt);
  }
}

\tikzdeclarepattern{
  name=primeddots2,
  type=uncolored,
  bounding box={(-.3pt,-.3pt) and (.3pt,.3pt)},
  tile size={(2.5pt,2.5pt)},
  tile transformation={rotate=60},
  parameters={none}, 
  code={
    \fill(0pt,0pt) circle (.35pt);
  }
}

\tikzstyle{primedregion}[none]=[
	preaction={fill=#1},
	pattern=primeddots,
  draw=#1,
]

\tikzstyle{primedregion2}[none]=[
	preaction={fill=#1},
	pattern=primeddots2,
  draw=#1,
]

\tikzcdset{scale cd/.style={every label/.append style={scale=#1},
    cells={nodes={scale=#1}}}}

\let\OLDthebibliography\thebibliography
\renewcommand\thebibliography[1]{
  \OLDthebibliography{#1}
  \setlength{\parskip}{0pt}
  \setlength{\itemsep}{0pt plus 0.3ex}
}

\begin{document}

\title{ {\bf Remarks on  \texorpdfstring{$\rm{C}^*$}{C*}-discrete inclusions}}
\author{Roberto Hern\'{a}ndez Palomares\footnote{Department of Mathematics, University of Waterloo \hfill \url{robertohp.math@gmail.com}}\hspace{.1cm} and Brent Nelson\footnote{Department of Mathematics, Michigan State University \hfill \url{brent@math.msu.edu}}}
\date{}
\maketitle
\begin{abstract}
\noindent We obtain a Galois correspondence between the lattice of intermediate C*-discrete subalgebras intermediate to a given irreducible C*-discrete inclusion, and characterize these as targets of compatible expectations under a traciality condition. Furthermore, we define \emph{freeness} for actions of tensor categories on C*-algebras, and show simplicity is preserved under taking reduced crossed products.
Finally, we show under certain constraints that $A$-valued semicircular systems give rise to C*-discrete inclusions, and thus are crossed products by an action of a tensor category. Along the way, we show the set of single algebraic generators of a dualizable bimodule forms an open subset. 
\end{abstract}

\section*{Introduction}
Irreducible C*-discrete inclusions were introduced in \cite[Theorem A]{2023arXiv230505072H}, and were characterized in terms of outer actions of unitary tensor categories (UTC s) over unital C*algebras with trivial center and a connected C*-algebra object.
As a consequence, these inclusions arise from an underlying C*-dynamical system governed by a quantum symmetry, and resemble the reduced crossed product inclusions by outer actions of discrete groups on unital C*-algebras with trivial center.  

Indeed, a unital inclusion $A\overset{E}{\subset}B$, where $E:B\twoheadrightarrow A$ is a faithful conditional expectation is irreducible (i.e. $A'\cap B\cong \bbC1$) and C*-discrete if and only if there is an \emph{outer action} $F:\cC\to \fgpBim(A)$ of $\cC$ on $A$, where $\cC$ is a unitary tensor category (UTC) and $F$ a fully-faithful unitary tensor functor, a $\cC$-graded C*-algebra object $\bbB$ which is \emph{connected} (ie $\bbB(1_\cC)\cong \bbC$),
and there exists an expectation preserving $*$-isomorphism 
$$
\left(A\overset{E}{\subset} B\right) \cong \left(A\overset{E'}{\subset} A\rtimes_{r}\bbB\right).
$$  
Here, $E'$ is the projection onto the $1_\cC$-graded component. (For details on $\cC$-graded algebra objects, see \cite{JP17, MR3948170}.)
As such, irreducible C*-discrete inclusions conform into a vast class of C*-inclusions determined by a standard invariant akin to subfactors.

The quantum symmetries govern the structure of an irreducible C*-discrete inclusion  through a Galois correspondence, which under natural assumptions parameterizes intermediate C*-discrete subalgebras $A\subseteq D\subseteq B$ that are ``compatible with the expectation'' by C*-algebra objects:
\begin{thmalpha}[{Corollary~\ref{cor:Galois}}]\label{thmalpha:Galois}
    Let $\cC$ be a unitary tensor category, $A$ be a unital C*-algebra with trivial center and an outer action $F:\cC\to\fgpBim(A).$
    Let $\left(A\overset{E}{\subset} B\right)\in\CDisc$ be an irreducible unital inclusion supported on $F[\cC]$, with corresponding connected $\cC$-graded C*-algebra object $\bbB\in\rCorr(\cC)$. 
    Then, there is a lattice injection: 
    \begin{align*}
        \left\{ D\in\rCorr\ \middle|\ 
            \begin{aligned}
                & A\subseteq D\subseteq B,\ \ \exists\text{ c.e. } E^B_D:B\twoheadrightarrow D,\\
                & \hspace{1.5cm}  E= E|_D\circ E^B_D
            \end{aligned}
        \right\}
         \hookrightarrow&
        \left\{\bbD\in\rCorr(\cC) \middle|\ \bbOne\subseteq \bbD\subseteq \bbB
        \right\}.
    \end{align*}
    Here, all C*-algebra inclusions are assumed unital. 
    In case $\bbB$ is tracial (see Equation \ref{eqn:TracialObject} for a definition), this injection is in fact a lattice isomorphism.
\end{thmalpha}
This statement is a state-free analogue of the Galois correspondence for discrete subfactors proven by C. Jones and Penneys in \cite[\S7.3]{MR3948170}. 
Simultaneously, it is a generalization of the Galois correspondence for an inclusion $A\overset{E}{\subset} A\rtimes_r\Gamma$ arising from the outer action of a discrete group $\Gamma$ on a unital \emph{simple} C*-algebra $A$ developed by Cameron and Smith in \cite[Theorem 3.5]{MR4010423}. 
Under these circumstances they showed  every intermediate C*-algebra is of the form $D=A\rtimes_{r}\Lambda$ for a subgroup $\Lambda\leq\Gamma$, and so is automatically $A$-discrete and compatible.  
In Section \ref{sec:GalDisc}, we show how to categorify some of their main results and incorporate them into our setting (cf Proposition \ref{prop:FullvsProjective}). 

We do not know yet if the hypotheses of discreteness or compatibility of expectations for $D$ are redundant in the setting of Theorem \ref{thmalpha:Galois} if simplicity of $A$ is assumed. 
In fact, Mukohara showed in \cite[Corollary 4.12]{2024arXiv240113989M} that for the fixed-points inclusion of an \emph{isometrically shift-absorbing} action of a compact group $G$ on a simple $A$, then every intermediate inclusion has the form $A^H$ for $H\subseteq G$ closed. 

\medskip

We now discuss when is a C*-discrete extension of a simple C*-algebra is simple.  
Since irreducible C*-discrete inclusions are crossed products by $\cC$, we sought for conditions in terms of the action $\cC\curvearrowright A.$
It was shown by Kishimoto \cite{MR634163} that if a discrete group $\Gamma$ acts outerly on a unital simple C*-algebra $A$, then the resulting reduced crossed product $A\rtimes_r \Gamma$ is also  simple. 
The results of \cite{MR4010423}, recapture Kishimoto's by establishing properties of certain convex neighborhoods of points in $A$ under the action of $\Gamma$. 
They further show these properties are enjoyed by all outer actions of discrete groups, and so we take it as our definition of \emph{freeness} for an action $\cC\curvearrowright A$ (Definition \ref{defn:FreeAction}).
In this article, we adapt Cameron and Smiths's simplicity result to \emph{free actions of UTCs}:
\begin{thmalpha}[Theorem \ref{thm:CrossProdSimple}] \label{thmalpha:CrossProdSimple}
    If $\cC \to\fgpBim(A)$ is a free action of a UTC on the simple unital C*-algebra $A$, and if $\bbB\in \rCorr(\cC)$ is any connected $\cC$-graded C*-algebra object, then the reduced crossed product
    $A\rtimes_{r}\bbB$
    is simple. 
\end{thmalpha}
\noindent The proof is analogous to Cameron and Smith's, where the freeness condition gives a firm control on the Fourier series for elements in the crossed product as in \cite[Proposition 3.3]{MR4010423}, allowing to export the simplicity of $A$ elsewhere. 

\medskip

In Section~\ref{sec:A-valued_semicirc} we investigate when $A$-valued semicircular systems, which arise frequently in free probability contexts, give irreducible C*-discrete inclusions. These C*-algebras $\hat{\Phi}(A,\eta)$ were introduced by Shlyakhtenko in \cite{Shl99}, and they are determined by a C*-algebra $A$ and a completely positive map $\eta$ from $A$ into an amplification of $A$. We provide sufficient conditions on $\eta$ to guarantee that the inclusion $A\subset \hat{\Phi}(A,\eta)$ is C*-discrete (see Theorem~\ref{thm:A-valued_discrete_inclusions}). It remains an open problem to fully determine the structure of the underlying dynamical system in terms of $\eta$. 

The above sufficient conditions rely on the stability of the algebraic generators of an $A$-$A$ bimodule under small perturbations, which is a problem of its own interest regarding C*-correspondences. 
We elaborate on this in Section~\ref{sec:Generators}, where in Theorem~\ref{thm:ProjectiveVectors} we establish that the set of single topological generators of an $A$-$A$ bimodule forms an open set. 
The structure of the set of single algebraic generators of $K\in \fgpBim(A)$ is closely related to that of $A$, and we expand on these aspects in Examples \ref{ex:uncomplementable}, \ref{ex:EndomorphicBimsStrictlyPositiveElts}, and Corollary \ref{cor:MeagerResidualVecors}.

\subsection*{Acknowledgements}
The authors are grateful to David Larson, Roger Smith, Reiji Tomatsu, and Matthew Kennedy who provided useful comments for the development of this project. We are indebted to Matthew Lorentz who provided tremendous support in the initial stages of this project. Special thanks to Miho Mukohara who noticed an error in the proof of Theorem A in a prior version of this manuscript and helped us correct it by highlighting the role of traciality for algebra objects. 
RHP was partially supported by the AMS-Simons Travel Grant 2022, NSF grant DMS-2001163, NSF grant DMS-2000331 and NSF grant DMS-1654159. BN was supported by NSF grant DMS-1856683.

\section{Preliminaries}\label{sec:Prelims}

The main abstract setting of this paper involves unital inclusions of C*-algebras 
$$
A\overset{E}{\subset}B,
$$
where $E:B\twoheadrightarrow A$ is a faithful conditional expectation. A canonical construction gives a $B$-$A$ right C*-correspondence $\cB$, which is the completion of $B$ as a right Hilbert $A$-module with $A$-valued inner product $\langle b_1|\ b_2\rangle_A:= E(b_1^*b_2).$ Automatically, there is a left $B$-action on $\cB_A$ by adjointable operators given by left multiplication. Moreover, since $E$ is assumed faithful, there is a faithful copy of $B$ inside of $\cB$, which we will systematically denote by $B\Omega\subset \cB,$ where $\Omega$ is the image of $1$ under this identification. And so we obtain $\overline{B\Omega}^{\|\cdot\|_A}= \cB.$

We shall now recall the definition the action of $\cC$, a UTC, on a unital C*-algebra $A$, and provide details of examples of interest. 
For a unital C*-algebra, $\fgpBim(A)$ denotes the UTC of finitely generated projective (fgp) $A$-$A$ bimodules with bimodular maps, and if $(N, \tau)$ is a $\rm{II}_1$-factor, $\spfgpBim(N)$ denotes the UTC of fgp $N$-$N$ \emph{spherical} bimodules with normal bimodular maps. The adjective spherical refers to those bimodules whose left and right \emph{von Neumann dimensions} match, and are sometimes referred to in the literature as \emph{extremal}. 
An action of  $\cC$, a UTC on $A$ (resp. on $N$) is given by a unitary tensor functor 
\begin{align*}
    &F:\cC\to\fgpBim(A) && (\text{respectively } F:\cC\to \spfgpBim(N)).
\end{align*}
We sometimes denote this by $\cC \overset{F}{\curvearrowright}A.$ 
For general properties, notation and expanded definitions for unitary tensor categories (UTCs), unitary tensor functors and $\cC$-graded $*$-algebra objects ($\Vec(\cC)$), we refer the reader to \cite{2023arXiv230505072H} and references therein.

Recall that a functor is said to be \emph{faithful} if it is injective at the level of hom spaces, and is said to be \emph{full} if it is surjective at the level of hom spaces.
By untarity, faithfulness of $F$ is automatic (since the relevant endomorphism spaces are finite dimensional C*-algebras, and $F(1)=1$). 
However, the property of being full is far from automatic,  and so we introduce the following definition:
\begin{defn}\label{defn:OuterAction}
    An \textbf{outer action} of a UTC $\cC$ on a unital C*/W*-algebra is a fully-faithful unitary tensor functor $F:\cC\to \fgpBim(A)$.
\end{defn}

This nomenclature has been used before in the literature, and in order to justify this nomenclature, we explain how to categorify an outer action by a discrete group on a unital C*-algebra:
\begin{ex}\label{ex:HilbgammOuterAction}
    For any countable discrete group $\Gamma$, we can consider the UTC of finite dimensional $\Gamma$-graded Hilbert spaces $\fdHilb(\Gamma),$ whose morphisms are uniformly bounded bounded linear transformations that respect the grading. 

    An action $\Gamma \overset{\alpha}{\curvearrowright} A$  over a simple unital C*-algebra $A$ consists of the same data as a  unitary tensor functor: 
    \begin{align*}
        \alpha:\fdHilb(\Gamma)&\to\fgpBim(A) \qquad \text{ with tensorator }  &\alpha^2_{g,h}:{}_{g}A\underset{A}{\boxtimes} {}_{h}A \to {}_{gh}A\\
        g&\mapsto {}_{g}A,  &{\xi\boxtimes\eta}\mapsto \alpha_{h^{-1}}(\xi)\cdot\eta.
    \end{align*}
    Here, for each $g\in\Gamma,$ ${}_{g}A\in\fgpBim(A)$ is $A$ as a vector space, with Hilbert C* $A$-$A$ bimodule structure given by $a\rhd\xi\lhd c=\alpha_{g^{-1}}(a)\cdot\xi\cdot c,$ and right and left inner products $\langle\xi\mid  \eta \rangle^{g}_A:=\xi^*\cdot \eta$ and ${}_A\!\langle\xi,\ \eta \rangle^{g}:=\alpha_{g^{1}}(\xi\cdot \eta^*)$ respectively.  
    
    In the usual language of group-actions, if $\{u_g\}_{g\in\Gamma}$ are the usual unitaries implementing the action, ie $\alpha_g(\cdot) = u_g(\cdot) u_g^*,$ outerness for $\alpha$ means that $g\neq h$ implies $u_g\neq u_h.$ 
    So, ${_g}A$ really is the same as $u_gA.$ 
    At the categorified level, this manifests as $g\neq h$ implies ${_g}A\not\cong {_h}A,$ so the categorified action maps distinct group elements to distinct irreducible bimodules.  
    And in greater generality this is precisely achieved by a fully-faithful functor, which will map non-isomorphic irreducible objects in $\cC$ into non-isomorphic simple objects in its target.  
\end{ex}
    
    We shall now discuss further symmetries of operator algebras which do not come from groups. 
    Popa established that every \emph{$\lambda$-lattice} arises as the standard invariant of some finite index subfactor \cite{MR1334479}. 
    Later on, together with Shlyakhtenko \cite{MR2051399} they showed this can be done using the $\rm{II}_1$-factor $\mathsf{L}\bbF_\infty$ associated to the free group in countably many generators, making it a universal receptacle for actions of UTCs \cite{BHP12}. 
    Subsequently, Guionnet, Jones, Shlyakhtenko and Walker provided a new diagrammatic proof of Popa's construction  via Jones' Planar algebras \cite{GJS10, JSW10}.     
    
\begin{ex}\label{ex:GJS}
    Using the \emph{GJS construction}, the first-named author and Hartglass show in \cite{MR4139893} that any (countably generated) UTC acts outerly on some unital separable simple exact C*-algebra $A$ with a unique trace $\tr$, whose $K$-theory depends only on $\cC$. 
    The GJS construction relies on diagrammatic and free-probabilistic techniques and yields an $A$-valued semicircular system C*-algebra, whose corners provide the bimodules used to represent $\cC.$ We refer the reader to \cite[\S5.1]{2023arXiv230505072H} for more details. 
\end{ex} 

Many other examples of actions of UTCs on C*-algebras are known, including \cite{MR1228532,MR1604162, MR4328058, MR4717816} and the references therein.

\section{Remarks on quasi-normalizers and singly generated bimodules}\label{sec:Generators}

Let $A$ be a unital C*-algebra and $K\in\fgpBim(A)$ be an fgp bimodule. 
As such, $K$ has left and right finite Pimsner-Popa bases. 
The goal of this section is to study the stability of bases under perturbation, as well as to explore how can $K$ be finitely or even singly generated as an $A$-$A$ bimodule. 
These questions are of fundamental nature and interesting in their own right, as they relate to the underlying Banach space structure of $K$, and will moreover be used later in Section \ref{sec:A-valued_semicirc} when we show how to obtain C*-discrete inclusions from $A$-valued semicircular systems. 
Towards the end of this section we will return to studying  inclusions of operator algebras, and revisit the quasi-normalizers involved in the definition of discreteness for inclusions. 

\subsection{Stability of Pimsner-Popa bases under small perturbations}

Given vectors $\{\xi_i\}_{i\in I}\subset K$ we write $$
A\rhd\{\xi_i\}_{i\in I}\lhd A:= \left\{\sum_{i\in I_0} a_i\rhd \xi_i\lhd a'_i\ \Bigg|\ a_i, a'_i\in A\right\}, 
$$
to denote their $A$-$A$ algebraic span.
Here, the sums rank over finite subsets $I_0\subset I$. 
We say $\{\xi_i\}_{i\in I}$ {\bf algebraically generates} $K$ if $ A\rhd \{\xi_i\}_{i\in I}\lhd A=K.$  
Similarly, $\{\xi_i\}_{i\in I}\subset K$ {\bf topologically generates} $K$ whenever $\overline{A\rhd\{\xi_i\}_{i\in I}\lhd A} =K.$
If $\xi\in K$, we let 
$$K_\xi:= \overline{A\rhd\{\xi\} \lhd A}$$
denote the singly topologically generated $A$-$A$ bimodule spanned by $\xi$. 
The two modes for spanning bimodules are in general distinct, and also $K_\xi$ need not be a complementable sub-bimodule of $K$ as we will see in Example \ref{ex:uncomplementable}. 

We are not aware of general necessary or sufficient conditions to guarantee $K$ has any {\bf $A$-$A$-cyclic vectors}; this is, single $A$-$A$ algebraic or topological generators. 
Nevertheless, in Theorem \ref{thm:ProjectiveVectors} we prove that single $A$-$A$-topologically-cyclic generators of $K$ conform into an open subset, and that they coincide with the single $A$-$A$-algebraic-cyclic vectors. 
In certain restricted cases of interest, in Corollary \ref{cor:MeagerResidualVecors} we obtained examples of bimodules where every element is $A$-$A$-cyclic up to a meager set. However, the general situation is still poorly understood. 

\smallskip

The following statement, \textit{Mostow's Lemma},  establishes that algebraic generators in a given Hilbert $A$-module are stable under small perturbations. 
Its proof is an elegant application of the fact that the set of epimorphisms between Banach spaces among bounded linear maps is open, and it can be found in \cite[Lemma 2.7.3]{MR2125398}. 
We reproduce the result here for the reader's convenience:
\begin{lem}\label{lem:BasisApproximation}
	Let $A$ be a unital C*-algebra. 
	Let $K\in \fgpMod(A)$ be a right $A$-module and let $\{u_i \}_{i=1}^{n}\subset K$ be an algebraic generating set. Then there is some $\delta > 0$ such that if $\{{u'}_i \}_{i=1}^{n}\subset K$ satisfy $||u_i - {u'}_i || < \delta, $ then $\{{u'}_i \}_{i=1}^{n}$ also algebraically generate $K$.
\end{lem}

We now use Mostow's Lemma to show stability of single generators of an fgp bimodule:   
\begin{thm}\label{thm:ProjectiveVectors}
    Let $A$ be a unital C*-algebra, $\cM$ a right C* $A$-$A$ correspoondence, and $\xi\in \cM$ be such that $\overline{A\rhd\xi \lhd A}=:\cM_\xi\in\fgpBim(A).$
    Then 
        $$\xi\in \big\{\eta\in K_\xi|\ A\rhd \eta\lhd A = \cM_\xi\big\}\subseteq \cM_\xi \text{ is open}.$$
    In particular $\xi$ generates $\cM_\xi$ algebraically.
\end{thm}
\begin{proof}
    Assume $\xi\neq 0,$ as otherwise the result holds trivially. 
    Let $\xi\neq 0$ be as in the statement. 
    Since $\cM_\xi\in \fgpBim(A)$ there are left and right Pimsner--Popa basis $\{v_j\}_{j=1}^m\subset \cM_\xi$ and $\{u_i\}_{i=1}^n\subset \cM_\xi$, respectively. 
    Thus, $\cM_\xi=\sum_{i=1}^n u_i\lhd A = \sum_{j=1}^m A\rhd v_j$  (algebraically). 
    But then Mostow's Lemma~\ref{lem:BasisApproximation} directly implies $\xi$ generates $\cM_\xi$ algebraically as an $A$-$A$ bimodule. Thus, for each 
    $i=1,\hdots, n$ and each $j=1,\hdots, m$ there exist finite subsets $\{a_{i,k},\ a'_{i,k}\}_{k\in I}$, 
    $\{c_{j,\ell},\ c'_{j,\ell}\}_{k\in J}\subset A$ such that 
    for each $i=1,\hdots, n$ and each $j=1,\hdots, m$ we have (finite sums)
    $$u_i=\sum_{k\in I} a_{i,k}\rhd\xi\lhd a'_{i,k},\quad \text{ and }\quad v_j=\sum_{\ell\in J} c_{j,\ell}\rhd\xi\lhd c'_{j,\ell}.$$ 
    Notice that not all these summands are zero for each expression. 
    Let $\delta_r>0$ be given by Mostow's Lemma for $\cM_\xi$ as generated as a right $A$-module by $\left\{a_{i,k}\rhd\xi\lhd a'_{i,k}\right\}_{i=1,\hdots, n}^{k\in I},$ and  
    similarly, let $\delta_l> 0$ be given by Mostow's Lemma for $\cM_\xi$ as generated as a left $A$-module by $\left\{c_{j,k}\rhd\xi\lhd c'_{j,k}\right\}_{j=1,\hdots, m}^{\ell\in J}.$ 
    
    Set $\delta:=\min\{\delta_l,\delta_r\}>0.$ 
    Now let $\xi'\in K_\xi$ be such that 
    $$0<||\xi-\xi'||_A < \delta\Bigg/\left( \sum_{i=1}^n\sum_{k\in I} ||a_{i,k}||\cdot||a'_{i,k}|| +  \sum_{j=1}^n\sum_{\ell\in J} ||c_{j,\ell}||\cdot||c'_{j,\ell}||\right).$$
    For each $i=1,\hdots, n$ we have the inequality
    $$\left\|u_i- \sum_{k\in I}a_{i,k}\rhd\xi'\lhd a'_{i,k}\right\|_A \leq \sum_{k\in I}||a_{i,k}||\cdot||a'_{i,k}||\cdot||\xi-\xi'||_A < \delta.$$
    Similarly, for each $j=1,\hdots, m$ we have 
    $$\left\|v_j - \sum_{\ell\in J}c_{j,\ell}\rhd\xi\lhd c'_{j,\ell}\right\|_A<\delta.$$ 
    Thus, further left and right application of Mostow's Lemma to the left/right generating subsets $\left\{c_{j,\ell}\rhd \xi'\lhd c'_{j,\ell}\right\}_{j=1,\hdots, m}^{\ell\in J}\subset \cM_\xi$ and  $\left\{a_{i,k}\rhd\xi \lhd a'_{i,k}\right\}_{i=1,\hdots, n}^{k\in I}\subset \cM_\xi$
    for our chosen $\delta$ implies each of them generates $\cM_\xi$ as a left and right $A$-module, respectively. This is, $A\rhd\xi'\lhd A=\cM_\xi,$ too.
\end{proof}
\noindent In a way, Theorem  \ref{thm:ProjectiveVectors} is a bimodular version of \emph{Kadison's Transitivity Theorem} applied to a cyclic Hilbert space representation.

Notice that if $\cM=A=K$ above, this recovers the well-known result for unital C*-algebras that  invertible elements $GL(A)\subset A$ conform into an open subset. 
Furthermore, Theorem \ref{thm:ProjectiveVectors} can be viewed as a generalized, or perhaps a higher version of (algebraic) simplicity for fgp $A$-$A$ bimodules. 
In Examples \ref{ex:uncomplementable} and \ref{ex:EndomorphicBimsStrictlyPositiveElts}, we shall illustrate cases of fgp $A$-$A$ bimodules where the set of single generators is non-trivial.

\begin{exs}[{\bf uncomplementable bimodules and $A$-$A$-cyclic vectors}]\label{ex:uncomplementable}
\ \\
    $\bullet$ Consider the unital C*-algebra $A=C[0,1],$ and the function $(\xi:t\mapsto t)\in A.$ We have that $ \xi A\subset \{f\in A|\ f(0)=0\}\subset A$, but clearly $\overline{\xi A}\neq A$ and $\overline{\xi A}^\perp =\{0\};$ i.e. $\overline{\xi A}$ is not complementable.  Furthermore, by the theory of bounded operators on Hilbert C*-modules (cf \cite{MR1325694, MR2125398},  \cite[Corollary 15.3.9]{wegge1993k}), $\overline{\xi A}$ is therefore not the range of any adjointable $A$-linear endomorphism. 
    Thus, $\overline{\xi A}$ is an example of a nontrivial and singly topologically generated Hilbert $A$-$A$ sub-bimodule of $A$ that is not fgp, as the latter are automatically complementable.

    On the contrary, if $\eta\in C[0,1]$ is invertible, then clearly $\eta A=A.$

    \noindent$\bullet \bullet$ Consider any unital non-simple C*-algebra $A$, and let $\{0\}\neq I\subsetneq A$ be a closed two-sided ideal, and consider $A$ as an $A$-$A$-bimodule.
    By construction, no vector in $I$ is $A$-$A$-cyclic for $A$. 
    One could consider the concrete case of the compact operators on the separable Hilbert space $I=\cK\subsetneq \End(\ell^2\bbN)=A$, where even when $Z(A)=\bbC\id$, no compact operator is $A$-$A$-cyclic for ${}_AA_A.$
    Notice $\cK^\perp=\{0\},$ and so $\cK$ is not a complementable summand inside $A.$ 
    It is moreover easy to see that any $k\in\cK^\times$ is a $A$-$A$-cyclic vector for ${}_{A}\cK_A$ 
\end{exs}

\begin{remark}\label{remark:fgisfgp}
     Every algebraically finitely generated Hilbert C*-module $\cM$ over a unital C*-algebra is fgp \cite[Thm 5.9]{MR1938798}; i.e. $\cM_A\cong p[A^{\oplus n}]_A,$ for some $n\in \bbN$ and a projection $p\in M_n(A).$
     As such, $\cM$ is always orthogonally complemented. 
     As we have seen, this feature is not shared by the finitely topologically generated Hilbert C*-modules. 
     For further remarks about these contrasts, as well as related questions, we direct the interested reader to \cite[\S 4]{2006math.....11349F}. 
\end{remark}

\subsection{Variants of quasi-normalizers for inclusions}

We now focus our attention to those bimodules arising from discrete inclusions. 
On the one hand, discreteness for a unital irreducible inclusion of C*-algebras $A\overset{E}{\subset} B$ was introduced in \cite{2023arXiv230505072H} by demanding that certain intermediate $*$-subalgebra $A\subset \PQN(A\subset B)\subset B$ called the \emph{projective quasi-normalizer} to be dense in $\|\cdot\|_B$-norm. 
As a shorthand, we use $B^\diamondsuit$ to denote $\PQN(A\subset B)$ and refer to it as {\bf the diamond space}. 
By definition, the projective quasi-normalizer of $A\subset B$ is 
\begin{equation}\label{eqn:diamondspace}
    B^\diamondsuit=\left\{b\in B|\ \exists K\in\fgpBim(A),\ b\Omega\in K\subset B\Omega \right\}. 
\end{equation}
On the other hand, \emph{quasi-regularity} for an \emph{irreducible extremal} subfactor $(N, \tau)\overset{E}{\subset} M$, where $N$ is type $\rm{II}_1$ with its unique trace $\tau$, as defined in \cite{MR3801484, MR3948170}, demands that $\QN(N\subset M)$, the \emph{quasi-normalizer} intermediate $*$-subalgebra $N\subset \QN(N\subset M)\subset M$ to be strongly dense in $M$ in the GNS representation $L^2(M, \tau\circ E)$.
To ease notation, we denote $\QN(N\subset M)$ by $M^\circ$. 
Here, extremality means that the fgp $N$-$N$ bimodules appearing in the semi-simple decomposition of the GNS space $L^2(M, \tau\circ E)=:L^2(M)$ as $N$-$N$ bimodules have matching left and right von Neumann dimensions. 
By definition, the quasi normalizer for $N\subset M$ is
\begin{equation}\label{eqn:QN}
    M^\circ=\left\{m\in M\Bigg|\ \exists \{x_i\}_1^R, \{y_j\}_1^{R'}\subset M,\ mN\subset \sum_1^R Nx_i,\ Nm\subset \sum_1^{R'}y_jN \right\}.
\end{equation}

As it turns out, the quasi-normalizer of $N\subset M$ is quite malleable and can be presented in various different ways. 
For example, by \cite[Lemmma 6.3]{FalguieresThesis}, $M^\circ$ can be re-phrased as the set of $m\in M$ such that 
$$
L_m:=\overline{\mathsf{span} \{\sum n_i\rhd m\Omega\lhd n'_i\}}^{\|\cdot\|_2}\in \spbfBim(N),
$$
where the sums rank over finite sets and the $n_i, n'_i\in N,$ $M\Omega$ is the faithful copy of $M$ inside $L^2(M),$ and the $\|\cdot\|_2$-norm is with respect to the faithful state $\tau\circ E$ on $M$. 
In other words, $M^{\circ}$ is the subset of $M$ whose elements generate an fgp $N$-$N$ bimodule.
Yet another presentation, following \cite[Lemma 2.5]{MR3801484} and \cite[Proposition 3.11]{MR3948170}, yields 
\begin{equation}\label{eqn:AlgRealization}
     M^\circ \cong \bigoplus_{\substack {L\in\Irr(\spfgpBim(N))\\ L\leq L^2(M)}}(L^{\oplus n_L})^\circ, 
\end{equation}
where $L^{\oplus n_L}$ is the isotypic component of $L$ in the semisimple decomposition of  ${}_NL^2(M)_N.$ 
Here, $(L^{\oplus n_L})^\circ :=(L^{\oplus n_L})\cap M\Omega$ is shown to correspond to the \emph{bounded vectors} in $L^{\oplus n_L}$, and thus the quasi-normalizer $M^\circ$ coincides with the space of all bounded vectors in $L^2(M).$

It was shown in \cite{MR3801484, MR3948170} that $N\subset M$ is quasi-regular if and only if it is discrete, under the assumption of irreducibility $N'\cap M= \bbC1$. 
In this context, discreteness means that the semi-simple decomposition into irreducible $N$-$N$ bimodules 
\begin{equation}\label{eqn:L2PW}
{}_NL^2(M)_N\cong \overline{\bigoplus_{L\in \Irr(\spfgpBim(N))} L^{\oplus n_L} }^{\|\cdot\|_2}   
\end{equation}
satisfies that every $n_L\in\bbN\cup\{0\},$ with the multiplicity of $L^2(N)$ be given by $n_N=1$. 
The equivalence of quasi-regularity and discreteness for irreducible subfactors allows one to work inside the Hilbert space $L^2(M)$ to conclude information about $M$.
This equivalence, along with the various presentations of $M^\circ$, are afforded by familiar Hilbert space facts such as bounded maps are automatically adjointable, and so closed subspaces are necessarily orthogonally complementable and are given by the range of self-adjoint idempotent endomorphisms (i.e. projections in $\cB(\cH)$). 
Of course, one is often forced to keep track of the \emph{spacial} data given by the faithful state $\tau\circ E$ on $M,$ which is why spherical bimodules appear in the first place, and which can also impose important restrictions in certain cases (cf \cite{MR2561199}). 

\medskip

We now return to the more general setting of $A\overset{E}{\subset} B$ as in the beginning of the section, where we moreover do not require the data of a state on $A$. 
In exchange for this broader point of view we have disposed of Hilbert space representations of $B$ and their toolkit, and instead have to work more intrinsically, using only the $B$-$A$ correspondence $\cB$ associated with $E.$ (See Section \ref{sec:Prelims} for definitions.)
This change demands revision of quasi-normalizers for $A\subset B$ in the spirit of Equations (\ref{eqn:QN}) and (\ref{eqn:AlgRealization}) in comparison with the projective quasi-normalizer introduced in Equation (\ref{eqn:diamondspace}), in lights of the discussion carried out in Section \ref{sec:Generators} around $A$-$A$ cyclic vectors and ways to generate a bimodule. 
We have that 
\begin{align}\label{eqn:AllQuasiNormalizers}
    \QN(A\subset B):= &B^\circ=\left\{b\in B\Bigg|\ \exists \{x_i\}_1^R, \{y_j\}_1^{R'}\subset B,\ bA\subset \sum_1^R Ax_i,\ Ab\subset \sum_1^{R'}y_jA \right\}\\
    &\rotatebox{90}{$\subseteq$}\nonumber\\
    &B^\diamondsuit\nonumber\\
    &\rotatebox{90}{$\subseteq$}\nonumber\\
    &B^{\diamondsuit\diamondsuit}:=\Big\{  b\in B \Big|\ \overline{A\rhd b\Omega \lhd A} \in\fgpBim(A) \Big\}.\nonumber
\end{align}
Notice, contrary to their behaviour with subfactors, these three sets need not, and in fact do not match in general. 
Indeed, in lights of Examples \ref{ex:uncomplementable} and \ref{ex:EndomorphicBimsStrictlyPositiveElts}, an $A$-$A$ bimodule need not be everywhere singly generated, and so $B^{\diamondsuit\diamondsuit}\neq B^{\diamondsuit}.$ (Algebraic or topological considerations in the definition of $B^{\diamondsuit\diamondsuit}$ will not alter this conclusion, given Theorem \ref{thm:ProjectiveVectors}.) Moreover, the $A$-$A$ bimodule topologically generated by elements $b\in B^\circ$ as a closed subspace of $\cB$ need not be the range of an adjointable idempotent in $\End({}_A\cB_A)$. 

\begin{ex}\label{ex:EndomorphicBimsStrictlyPositiveElts}
    In the very specific case where an $A$-$A$ fgp bimodule comes from an endomorphism, we can ensure there indeed exist $A$-$A$-cyclic vectors. 
    In doing so, this will also provide concrete examples where the sets from Expression (\ref{eqn:AllQuasiNormalizers}) differ. 
    Extrapolating from previous notation, if $K\in \fgpBim$, we write $K^{\diamondsuit\diamondsuit}$ to denote the $A$-$A$-cyclic vectors $\xi\in K$; this is $\overline{A\rhd \xi \lhd A}\in \fgpBim(A)$.

    Recalling the notation from Example \ref{ex:HilbgammOuterAction}, consider the fgp $A$-$A$ bimodule ${}_\lambda A_A,$ where $\lambda:A\to A$ is a unital $*$-homomorphism and $a\rhd \xi\lhd a'= \lambda(a)\xi a'$ for $a,a',\xi\in A.$ Then, for an arbitrary $\xi\in A,$ the algebraic $A$-$A$ bimodule it generates is exactly the right ideal $\lambda[A]\xi A$. 
    If $\lambda \in\Aut(A),$ and $A$ is simple, then any $\xi$ will generate the entire fgp bimodule ${}_\lambda A_A$ algebraically. 
    In the case where only $\lambda\in\End(A)$ has finite index, since $\xi A\subseteq \lambda[A]\xi A$ by unitality, it is clear that any right-invertible $\xi$ will generate all ${}_\lambda A_A$ algebraically. 
    So at least 
    $$ GL(A)\subseteq {}_\lambda A_A^{\diamondsuit\diamondsuit}\neq\emptyset.
    $$ 
    If moreover, $\lambda\in \End(A)$ and $A$ is assumed simple, then we also have 
    $$
    \lambda[A]\subset {}_\lambda A_A^{\diamondsuit\diamondsuit}.
    $$ 

    Whenever $0\neq\xi\in {_\lambda}A_A\in\fgpBim(A)$, there are certain abstract criteria that would make $\overline{\lambda[A]\xi A}$ into a complementable subbimodule of ${_\lambda}A_A$ and thus equal to it. For example, if $\overline{\lambda[A]\xi A}$ is \emph{self-dual} it is necessarily orthogonally complementable \cite[Exercise 15.I(d)]{wegge1993k}. Also, if $\overline{\lambda[A]\xi A}$ contains a \emph{strictly positive element}---or $\xi$ itself is one---it then topologically generates $A$ as right ideal \cite[Chapter 6]{MR1325694}. (See references for definitions.) These conditions, however, might be hard to check in practice, but in either case, algebraic generation is sufficient by Theorem  \ref{thm:ProjectiveVectors}. 
\end{ex}

The better behaved among the three is $B^\diamondsuit$. 
By the C*-Frobenius Reciprocity Theorem \cite[Theorem 2.6]{2023arXiv230505072H}, it can be shown that $B^\diamondsuit$ is an algebraic realization of a C*-algebra object in $\Vec(\cC)$. Moreover, if the inclusion is irreducible $A'\cap B=\bbC1$, then $B^\diamondsuit$ admits a decomposition
\begin{equation}\label{eqn:DiamondDecomposition}
    B^\diamondsuit\cong \bigoplus_{\substack {K\in\Irr(\fgpBim(A))\\ K\leq \cB\\ K\subset B\Omega}}K^{\oplus n_K}, \qquad n_K\in \bbN\cup\{0\}
\end{equation}
and thus, $A\subseteq B^\diamondsuit\subset B$ is a *-algebra.
In comparison, since determining the $A$-$A$-cyclic vectors in $B\Omega\subset\cB$ that generate an fgp bimodule is currently out of reach, it is not obvious or known whether $B^{\diamondsuit\diamondsuit}$ is even closed under sums. 
Regarding $B^\circ$, the problem is that we do not have means to discriminate which $b\in B^\circ$ will generate a complementable $A$-$A$ sub-bimodule of $\cB$, and cannot make further claims about the algebraic structure of $B^\circ.$ 
In lights of Equation (\ref{eqn:DiamondDecomposition}) in comparison with Equations (\ref{eqn:AlgRealization}) and (\ref{eqn:L2PW}), one could interpret the $*$-algebra $B^\diamondsuit$ as a state-less replacement of bounded vectors.

Arguably, the choice of $B^\diamondsuit$ in the definition of discreteness is the correct one, as demonstrated by \cite[Theorem A]{2023arXiv230505072H}, establishing a correspondence between irreducible C*-discrete inclusions and actions of unitary tensor categories along with their connected C*-algebra objects, affording a meaningful version of standard invariants \cite[Corollary C]{2023arXiv230505072H}.
It seems likely that the open problem of explicitly characterizing the sets $B^{\diamondsuit\diamondsuit}$ and $B^\circ$ lies at the heart of the Galois correspondence studied in Section \ref{sec:Galois}, which uses the *-algebra $B^\diamondsuit$. And, in certain cases (cf \cite{MR4599249, 2024arXiv240113989M}), it might be possible to establish stronger conclusions, shedding light on the structures of $B^\circ$ and $B^{\diamondsuit\diamondsuit}$.

If one assumes that $A$ has certain internal structure, then it is possible to make stronger conclusions about its bimodules: 
\begin{cor}\label{cor:MeagerResidualVecors}
    In the setting of Example \ref{ex:EndomorphicBimsStrictlyPositiveElts}, if $A$ has \emph{topological stable rank one} and ${}_\lambda A_A\in\fgpBim(A),$
    then 
    $$
        {}_\lambda A^{\diamondsuit\diamondsuit}_A\subset A \text{ is dense}.    
    $$
    In particular, its complement is a meager set. 
    
\end{cor}
\begin{proof}
    Since $A$ has topological stable rank one, then by \cite[V.3]{MR2188261} then $GL(A)\subset A$ is dense, and so algebraic generators of ${}_\lambda A_A$ automatically form a dense subset.
    By Theorem \ref{thm:ProjectiveVectors}, it is also open, and so its complement is meager. 
\end{proof}
Remarkably, the GJS C*-algebras introduced in Example \ref{ex:GJS} have topological stable rank one \cite{MR3624399}. 

The bimodules appearing in Corollary \ref{cor:MeagerResidualVecors} typically arise from the study of \emph{superselection sectors}, which we think of as a unitary tensor category (assuming $Z(A)=\bbC1$) 
$\mathsf{Sect}_\mathsf{f}(A),$ whose objects are finite-index endomorphic bimodules ${}_\lambda A$ as above. 
Given finite-index endomorphisms $\rho,\sigma$ of $A$ their tensor product is given by reversed composition, and their intertwiners $\{a\in A| a\rho(\xi)= \sigma(\xi)a \forall \xi\in A\}$.

\section{Galois Correspondence}\label{sec:Galois} 
The main purpose of this section is to parameterize the lattice of certain intermediate inclusions to a given fixed irreducible $A\overset{E}{\subset}B\in\CDisc$ in terms of certain C*-algebra objects graded by the support of this inclusion, and in doing so establish a Galois correspondence in the spirit of \cite{MR1622812, MR1900138, MR4010423, MR3948170, MR4083877}. 

The main difficulty we encounter is the problem of determining whether $\CDisc$ is closed under taking unital subalgebras. This is, if $A\overset{E}{\subset}B\in\CDisc$ is it always true that $A\overset{E|_D}{\subset}D\in\CDisc$? 
We managed to partially answer this question positively in Theorem \ref{thm:ceCDisc} provided the associated C*-algebra object is tracial (see Equation \ref{eqn:TracialObject}) in combination with Theorem \ref{thm:ConverseceCDisc}, where $D$ is the target of a conditional expectation $E^B_D:B\twoheadrightarrow D$ that is compatible with $E$.

Examples of C*-discrete inclusions that are closed under taking subalgebras include all  finite-index extensions of simple unital C*-algebras, where in \cite[Proposition 6.1]{MR1900138} Izumi proves the existence of conditional expectation onto intermediate C*-subalgebras. Furthermore, in the context of outer actions of countable and discrete groups on unital simple C*-algebras, Cameron-Smith show in \cite[Theorem 3.5]{MR4010423} that the intermediate C*-algebras are always the range of some conditional expectation. However, in this latter case the fact that the bimodules supporting the inclusion are given by $*$-automorphisms plays a crucial role in the development of their results. 

We do not currently have a complete answer to the general problem, and, in fact, we have multiple reasons to believe there might exist certain \emph{residual} intermediate subalgebras which are not C*-discrete. 
For example, in \cite[Theorem 5.1]{MR4083877}, Suzuki shows that purely infinite algebras always have a infinite index subalgebra given by the image of an \emph{exotic endomorphism}, which is not the range of any conditional expectation. Another reason to believe that not all intermediate C*-algebras $A\subseteq D\subseteq B$ are targets of compatible conditional expectations follows from the fundamental problem of determining when an irreducible $K\in\fgpBim(A)$ is everywhere singly algebraically generated (see Section \ref{sec:Generators}) provided $A$ is simple. We shall defer further exploration of these problems for later work, and restrict ourselves to study the targets of conditional expectations. Several instances of this behaviour where highlighted in Examples \ref{ex:uncomplementable} and \ref{ex:EndomorphicBimsStrictlyPositiveElts}. 
Interestingly, the fgp bimodules supporting the cores of Cuntz algebras are algebraically generated by any nonzero element  \cite[Proposition 3.4 (4)]{2025arXiv250321515H}. 

\begin{remark}
    In the context of subfactors arising from minimal actions of compact groups on infinite factors, Izumi-Longo-Popa show in \cite[Theorem 3.9]{MR1622812} that intermediate subfactors are always targets of some conditional expectation.
    Tomatsu further generalized this setting to include minimal actions of compact quantum groups, and established in \cite[Theorem 1]{MR2561199} that the discrete subfactors arising this way are closed under taking intermediate inclusions. 
    We will further comment on Tomatsu's results in Remark \ref{remark:TomatsuHaarPodles}.

    Moreover, irreducible projective bimodules over a $\rm{II}_1$-factor are always singly algebraically generated and the fundamental problem of having enough single bimodule generating vectors does not appear. This is a manifestation that the representation categories of factors are precisely unitary tensor categories, while the representation categories for trivial-center C*-algebras might be substantially more complex. Take for example ultra-weakly closed one-sided ideals in a factor, which are always principal and generated by a projection \cite[Proposition 3.12]{MR1873025}. In contrast, a principal one-sided closed ideal in a C*-algebra might often fail to be generated by a projection! (See for instance \cite[Exercise 15L]{wegge1993k} for a spectral criteria.)  
\end{remark}

    If $(\bbB, j^\bbB)$ is a connected C*-algebra object, the vector spaces $\bbB(a)$ for $a\in \Irr(\cC)$ come equipped with nondegenerate left and right inner products: for $f,g\in \bbB(a)$ define 
    \begin{align}\label{eqn:InnerProducts}
    \begin{split}
        {}_a\langle f,\ g\rangle &:= \bbE_{1_\cC}(\bbB^2_{a, \bar{a}}(f\otimes j^\bbB(g)))\\
        \langle g\ |\ f\rangle_a &:= \bbE_{1_\cC}(\bbB^2_{\bar{a},a}(j^\bbB(g)\otimes f)),
    \end{split}
    \end{align}
    where $\bbE_{1_\cC}$ denotes the projection onto the ground algebra 
    $\bbB(1_\cC)\cong \bbC.$ 
    Thus, for every $a\in \Irr(\cC),$ $\bbB(a)$ has two natural Hilbert space structures, which are moreover finite-dimensional by \cite[Proposition 2.8]{MR3948170}.

    We say that $\bbB$ {\bf is tracial} if and only if 
    \begin{align}\label{eqn:TracialObject}
    {}_a\langle f,\ g\rangle = \langle g\ |\ f\rangle_a
    \end{align}
    for every $a\in \Irr(\cC)$ and each $f,g\in \bbB(a)$  \cite[Definition 4.32]{JP17}.   
    Notice that if $\bbD\subseteq \bbB$ is a subalgebra of a tracial C*-algebra object, then $\bbD$ is automatically tracial, and we can write 
    $$
    \bbB(a)\cong \bbD(a)\oplus \bbD^\perp(a),
    $$
    where the orthogonal complement of $\bbD(a)$ inside $\bbB(a)$ does not depend on the choice of inner product.
\noindent We say {\bf an irreducible C*-discrete inclusion $A\subset B$ is tracial} if and only if the associated connected C*-algebra object $\bbB$ is tracial. 
%
%
    This naturally lifts to an orthogonal projection at the algebra object level (i.e. a ucp categorical map)
    $$P^\bbB_\bbD:\bbB\Rightarrow\bbB,$$ with $P^\bbB_\bbD[\bbB]=\bbD$ and thus 
    \begin{align}\label{eqn:bbBDecomp}
        \bbB=\bbD\oplus\bbD^\perp.
    \end{align}

The following theorem establishes that the intermediate discrete subalgebras to an irreducible C*-discrete tracial inclusion are automatically targets of compatible conditional expectations: 
\begin{thm}\label{thm:ceCDisc}
    Let $\cC$ be a unitary tensor category, $A$ be a unital C*-algebra with trivial center and an outer action $F:\cC\to\fgpBim(A).$
    Let $\left(A\overset{E}{\subset} B\right)\in\CDisc$ be an irreducible unital inclusion supported on $F[\cC]$. 
    Let $\bbB\in\rCorr(\cC)$ be the associated connected C*-algebra object from \cite[Theorem A]{{2023arXiv230505072H}} (i.e. $B\cong A\rtimes_{r,F}\bbB$), and let $\bbD\in\rCorr(\cC)$ with $\bbOne\subseteq\bbD\subseteq \bbB$ and $D:=A\rtimes_{r,F}\bbD$. 
    If $\bbB$ is tracial, then there exists a faithful conditional expectation 
    \begin{align*}
        E^B_D:B\twoheadrightarrow D
    \end{align*}
    satisfying $E= E|_D\circ E^B_D$.
\end{thm}
\begin{proof}
%
    Consider the orthogonal decomposition from Equation \ref{eqn:bbBDecomp}, where the choice of inner product on the fibers of $\bbB$ is immaterial by traciality.
    The \emph{algebraic realization isomorphism} from \cite[Proposition 2.5]{{2023arXiv230505072H}}  yields the orthogonal direct sum decomposition as $A$-$A$ right C*-correspondences 
    $$\cB = \overline{\PQN(A\subset B)\Omega }^{\|\cdot\|_A } \cong (A\rtimes_{F}\bbD)\oplus \left(A\rtimes_{F}\bbD^\perp\right) = \cD\oplus \cD^\perp.$$ 
    Let $p\in\End^\dag({_A}\cB_A)$ be the orthogonal projection such that $\cD\cong p[\cB]$ and $(1-p)[\cB]\cong \cD^\perp.$ 
    In light of the algebraic isomorphism from \cite[Proposition 2.5]{{2023arXiv230505072H}}, it is evident that $p$ matches with $\id_F\otimes P^\bbB_\bbD$ from the discussion preceeding Equation \ref{eqn:bbBDecomp},
     when restricted to $B^{\diamondsuit}\Omega.$  
    And in fact, if we denote $D^{\diamondsuit}:=\PQN\left(A\overset{E|_D}{\subseteq}D\right)$, by the C*-discreteness hypothesis $B=\overline{\PQN(A\subset B)}^{\|\cdot\|_B}$ it then follows that $\overline{\check{p}[B^{\diamondsuit}\Omega]}^{\|\cdot\|_B} = D$,  
    since $p[\PQN(A\subset B)\Omega]=p[B^{\diamondsuit}\Omega] = D^{\diamondsuit}\Omega$.
    (Recall that $\check{\ }$ deletes $\Omega$ as in \cite[Definition 2.4]{2023arXiv230505072H}.)

    We are now ready to construct $E^B_D.$ For each $b\in B^\diamondsuit$, by the faithfulness of $E$, there is a unique $E^B_D(b)\in D^{\diamondsuit}$ satisfying $p(b\Omega) = E^B_D(b)\Omega$ by the previous conclusions. 
    Now, recall that there is a faithful unital $*$-homomorphism $-\rhd-:B\to \End^\dag(\cB_A),$ given by (the extension of) multiplication on the left by elements of $B$, and we identify $B$ and its sub-algebra $D$ with their images. 
    We shall now prove that
    \begin{align}\label{eqn:Implementation}
        E^B_D(b)pb'\Omega = pbpb'\Omega \qquad \forall b,b'\in B^\diamondsuit    
    \end{align}
    so that $p$ implements $E^B_D$ on $\cB$.

    Denote by $(D^\diamondsuit)^\perp:= (1-p)\check{}\ [B^\diamondsuit\Omega]$. We observe that $j^\bbB[\bbD(c)^\perp]= \bbD(\overline{c})^\perp$ since $j^\bbB[\bbD(c)]=\bbD(c)$ and $j^\bbB$ is an involution.
    Now we make the following assertion:\\
    \underline{Claim 1}: $(D^\diamondsuit)^\perp$ is $D^\diamondsuit$- $D^\diamondsuit$ invariant.\\
    \underline{Proof}: 
    For elementary tensors $\xi_a\otimes f_a,\ \zeta_c\otimes h_c\in D^\diamondsuit$ and $\eta_b\otimes g_b^\perp\in (D^\diamondsuit)^\perp$ for $a,b,c\in \Irr(\cC)$ we have: 
    \begin{align*}
        \langle \xi_a\otimes f_a\Omega\ |\ (\eta_b\otimes g_b^\perp)\cdot(\zeta_c\otimes h_c)\Omega\rangle^\cB_A &= \left(\langle f_a\ |\ \ g_b^\perp\cdot h_c\rangle^\bbB \right)   \cdot \left( \langle \xi_a\ |\ \eta_b\cdot \zeta_c\rangle^{\bbF}_A\right)\\
        &= \bbE^\bbB_{1_\cC}\left(j^\bbB(f_a)\cdot (g_b^\perp\cdot h_c)\right)\cdot \bbE^\bbF_{1_\cC}\left(j^\bbF(\xi_a)\cdot (\eta_b\cdot \zeta_c)\right)\\
        &\overset{(\spadesuit)}{=} \bbE^\bbB_{1_\cC}\left(g_b^\perp \cdot (h_c\cdot j^\bbB(f_a))\right)\cdot \bbE^\bbF_{1_\cC}\left(j^\bbF(\xi_a)\cdot (\eta_b\cdot \zeta_c)\right)\\
        &=\left\langle j^\bbB(g_b^\perp)\ \big|\ h_c\cdot j^\bbB(f_a)\right\rangle\cdot \bbE^\bbF_{1_\cC}\left(j^\bbF(\xi_a)\cdot (\eta_b\cdot \zeta_c)\right) = 0.
    \end{align*}
    Notice that the $(\spadesuit)$ equation above holds by the traciality assumption on $\bbB.$ Moreover, the last equality follows since $(h_c\cdot j^\bbB(f_a))\in \bbD(c\otimes \bar{a})$ since $\bbD$ is a $*$-algebra object, and $j^\bbB(g_b^\perp)\in \bbD(\bar{b})^\perp$. Consequently, $(D^\diamondsuit)^\perp$ is right $D^\diamondsuit$-invariant by extending to finite sums. 

    To check left $D^\diamondsuit$-invariance we have
    \begin{align*}
        \langle \xi_a\otimes f_a\Omega\ |\ (\zeta_c\otimes h_c)\cdot(\eta_b\otimes g_b^\perp)\Omega\rangle^\cB_A &= \left(\langle f_a\ |\ \ h_c\cdot g_b^\perp \rangle^\bbB \right)   \cdot \left( \langle \xi_a\ |\ \zeta_c\cdot \eta_b\rangle^{\bbF}_A\right)\\
        &=\bbE^\bbB_{1_\cC}\left(j^\bbB(f_a)\cdot (h_c\cdot g_b^\perp)\right)\cdot \bbE^\bbF_{1_\cC}\left(j^\bbF(\xi_a)\cdot (\zeta_c\cdot \eta_b)\right)\\
        &=\bbE^\bbB_{1_\cC}\left((j^\bbB(f_a)\cdot h_c) \cdot g_b^\perp\right)\cdot \bbE^\bbF_{1_\cC}\left((j^\bbF(\xi_a)\cdot \zeta_c)\cdot \eta_b\right)\\
        &=\left\langle j^\bbB(j^\bbB(h_c)\cdot f_a)\ \big|\ g_b^\perp\right\rangle\cdot \bbE^\bbF_{1_\cC}\left((j^\bbF(\xi_a)\cdot \zeta_c)\cdot \eta_b\right) =0.
    \end{align*}
    Here, we solely relied on associativity of products and the closure properties of $ \bbD$ and its complement used above. This part argument works without traciality (and is moreover not strictly needed to prove the theorem).  
    We then have that $(D^\diamondsuit)^\perp$ is left $D^\diamondsuit$-invariant by extending to finite sums.

    \underline{Claim 2}: As on operators on $\cB$, $pdp=dp$ for $d\in D^\diamondsuit$ and $pdp=0$ for $e\in (D^\diamondsuit)^\perp$.\\
    \underline{Proof:}
    For $b\in B^\diamondsuit$ we have
        \[ 
            pdp b \Omega = p d E^B_D(b) \Omega = d E^B_D(b) \Omega = dpb\Omega,
        \]
    where we have used that $D^\diamondsuit$ is an algebra, and
        \[
            pepb\Omega = p e E^B_D(b)\Omega = 0,
        \]
    where we have used Claim 1 to assert that $e E^B_D(b)\in (D^\diamondsuit)^\perp$. Claim 2 then follows by extending to $\cB$ by continuity.

    We now prove Equation \ref{eqn:Implementation}. 
    For each $c\in \Irr(\cC),$ we let $\{f^{i_c}\}_{1\leq i_c\leq \dim(\bbB(c))}\subset \bbB(c)$ be an orthonormal basis, and $1\leq m_c\leq \dim(\bbB(c))$ be such that $\{f^{i_c}\}_{1\leq i_c\leq m_c}$ is an ONB for $\bbD(c).$ 
     Let $b = \sum_{c\in\Irr(\cC)}\sum_{i_c}\xi^{i_c}\otimes f^{i_c}\in B^\diamondsuit$ be arbitrary. 
     We have that 
    \begin{align*}
        pbp &= \sum_c\sum_{i_c=1}^{\dim(\bbB(c))} p(\xi^{i_c}\otimes f^{i_c})p\\
        &\overset{\clubsuit}{=} \sum_c\sum_{1\leq i_c\leq m_c} (\xi^{i_c}\otimes f^{i_c})p\\
        &= E^B_D(b)p.
    \end{align*}
    Here, the equation signaled by $\clubsuit$ follows from Claim 2.

    By Equation \ref{eqn:Implementation} we can therefore extend $E^B_D$ to a ucp map on $B=\overline{B^\diamondsuit}^{\|\cdot\|_B}$, which is $D$-$D$ bimodular. Using $\Omega\in \cD$, we then have
        \[
            E(b) = \< \Omega \mid b\Omega\>_A = \< \Omega \mid p b p \Omega\>_A = \<\Omega\mid E_D^B(b) \Omega\> = E|_D \circ E_D^B(b),
        \]
    for all $b\in B$. Thus the faithfulness of $E^B_D$ follows from that of $E$.
\end{proof}

We now establish a partial converse, which forces intermediate compatible subalgebras to an irreducible (not necessarily tracial) C*-discrete inclusion to be automatically C*-discrete. This gives an infinite index version of \cite[Lemma 6.2]{MR1900138}. 
\begin{thm}\label{thm:ConverseceCDisc}
    Let $\cC$ be a unitary tensor category, $A$ be a unital C*-algebra with trivial center and an outer action (i.e. a fully-faithful unitary tensor functor) $F:\cC\to\fgpBim(A).$
    Let $\left(A\overset{E}{\subset} B\right)\in\CDisc$ be an irreducible unital inclusion supported on $F[\cC]$. 
    If $A\subseteq D\subseteq B$ is an intermediate unital C*-algebra admitting a conditional expectation $E^B_D:B\twoheadrightarrow D$ that satisfies $E=E|_D\circ E^B_D,$ then there is a connected C*-algebra object $\bbD\in\rCorr(\cC)$ with $\bbOne\subseteq \bbD\subseteq\bbB$ and such that $D\cong A\rtimes_{F,r}\bbD$ so that $\left(A\overset{E|_D}{\subseteq} D\right)\in\CDisc.$ 
\end{thm}
\begin{proof}
    With the definitions as in the statement, for each $c\in\cC,$ define  
    \begin{align*}
        \bbD(c) &:= \left\{(E^B_D\circ \check{f}(\,\cdot\,))\Omega \ \Big|\ f\in \rCorr_{A-A}(F(c)\to \cB)^\diamondsuit=\bbB(c)\right\}\\
        &=\left\{g\in \rCorr_{A-A}(F(c)\to \cB)^\diamondsuit \ \Big|\ \ g[F(c)]\subset D\Omega\right\}\\
        &\cong \rCorr_{A-D}\left(F(c)\boxtimes_A D\to D\right).    
    \end{align*}
    The last isomorphism is an application of C*-Frobenius Reciprocity Theorem \cite[Theorem 2.6]{2023arXiv230505072H} to the unital inclusion $\left(A\overset{E|_D}{\subseteq} D\right)\!.$
    We claim this shows that $\bbD=\langle D\rangle$ is the underlying C*-algebra object to $A\subset D.$ To see this, it is clear that for every $c\in\cC$ we have $\langle D \rangle(c) \subset \bbD(c),$ since $E^B_D$ is the identity on $D$. The reversed containment follows from the assumption of compatibility: it is clear that each $E^B_D(\check{f}(\,\cdot\,))\Omega\in\bbD(c)$ is $A$-$A$ bimodular, and adjointability follows from
    $$\langle d\Omega\mid E^B_D(\check{f}(\xi))\Omega \rangle^\cB_A = E(d^*E^B_D(\check{f}(\xi)))=E(E^B_D(d^*\check{f}(\xi))) = E(d^*\check{f}(\xi)) = \langle f^\dag(d\Omega)\mid \xi\rangle^{F(c)}_A.$$
    Consequently, $\bbD\in\rCorr(\cC)$ is the connected C*-algebra object $\langle D\rangle,$ with $\bbOne\subseteq \bbD\subseteq \bbB$. 
    
    We now show that $D=A\rtimes_{F,r}\bbD$ and consequently that the inclusion $A\subset D$ is C*-discrete. By construction and \cite[Proposition 2.5]{2023arXiv230505072H}, $(A\rtimes_{F,r}\bbD)^\diamondsuit := \PQN(A\subset A\rtimes_{F,r}\bbD) \subset D,$ and so $(A\rtimes_{F,r}\bbD)^\diamondsuit \subseteq D^{\diamondsuit}:=\PQN(A\subset D)$.  
    Now, let $b\in B^{\diamondsuit}\setminus (A\rtimes_{F,r}\bbD)^\diamondsuit,$ where as usual $B^{\diamondsuit}:=\PQN(A\subset B).$ 
    Thus, for some minimal $L\in\fgpBim(A)$ we have that $b\Omega\in L\subset B\Omega.$ We now orthogonally decompose into irreducible fgp $A$-$A$ bimodules, $L\cong \bigoplus_1^m L_i^{\oplus n_i}$ where without loss of generality for some $i,$ some copy of $L_i$ is orthogonal to $(A\rtimes_{F,r}\bbD)^\diamondsuit\Omega$ (see the proof of \cite[Theorem 2.7]{2023arXiv230505072H}. 
    But then $L_i\perp D^{\diamondsuit}\Omega$ since $\bbD =\langle D\rangle,$ and so $b\not\in D^{\diamondsuit}.$ 
    Thus $D^{\diamondsuit}=(A\rtimes_{F,r}\bbD)^\diamondsuit,$ concluding the proof. 
\end{proof}

The combination of Theorems \ref{thm:ceCDisc} and \ref{thm:ConverseceCDisc} precisely characterizes C*-discrete subalgebras intermediate to a connected tracial C*-discrete inclusion as ranges of compatible expectations.
We turn our attention to parameterize this lattice in terms of C*-algebra objects and their morphisms in the spirit of \cite[Theorem A]{2023arXiv230505072H} involving the following categories: 
    \begin{itemize}
        \item The subcategory of $\CDisc$ whose morphisms consist of unital $A$-$A$ bimodular C*-algebra homomorphisms preserving the expectation, and
        \item The subcategory of $\cCCAlgs$ whose morphisms are natural transformations of $*$-algebra objects.
    \end{itemize}
    
We record the obtained Galois correspondence in the following corollary.
\begin{cor}[{Theorem~\ref{thmalpha:Galois}}]\label{cor:Galois}
    Let $\cC$ be a unitary tensor category, $A$ be a unital C*-algebra with trivial center and an outer action $F:\cC\to\fgpBim(A).$
    Let $\left(A\overset{E}{\subset} B\right)\in\CDisc$ be an irreducible unital inclusion supported on $F[\cC]$, with corresponding connected $\bbB\in\rCorr(\cC)$. 
    We then have the following lattice injection: 
    \begin{align*}
        \left\{ D\in\rCorr\ \middle|\ 
            \begin{aligned}
                & A\subseteq D\subseteq B,\ \ \exists\text{  c.e. }E^B_D:B\twoheadrightarrow D,\\
                &\hspace{1.5cm} E= E|_D\circ E^B_D
            \end{aligned}
        \right\}
         \hookrightarrow&
        \left\{\bbD\in\rCorr(\cC) \middle|\ \bbOne\subseteq \bbD\subseteq \bbB
        \right\}\\
         \cong &\ \left\{ D\in\rCorr\ \middle|\ 
            \begin{aligned}
                & A\subseteq D\subseteq B,\\
                & A\overset{E|_D}{\subset}D \in\CDisc
            \end{aligned}
            \right\}.
    \end{align*}
    Here, all C*-algebra inclusions are assumed unital. 
    In case $\bbB$ is tracial, then the injection is surjective, giving a lattice isomorphism. 
\end{cor}

\begin{remark}\label{remark:TomatsuHaarPodles}
The most natural examples of tracial C*-algebra objects come from the fiber functors associated to compact quantum groups of Kac-type \cite[Proposition 6.4]{MR3948170} (see Examples~\ref{ex:non/tracial_quantum_groups} below).
In an earlier version of this manuscript, we overlooked the role of traciality in the proof of Theorem \ref{thm:ceCDisc}, which we now have  clarified. 
In the tracial/Kac setting, there is therefore no intersection with Tomatsu's \cite[Example 3.13]{MR2561199} involving compact quantum groups of non Kac-type . 

\end{remark}

\begin{exs}\label{ex:non/tracial_quantum_groups}[{\bf Tracial and non-tracial connected C*-algebra objects}]
    Since the hypothesis of traciality plays a central role in Theorem \ref{thmalpha:Galois}, we highlight some of the examples of tracial connected C*-algebra objects—as well as non-examples—as developed in \cite{MR3948170}, where the relevant details can be found.
    
    Consider a discrete quantum group $\bbG$ as the pair $(\cC, F),$ where $\cC$ is a unitary tensor category, and $F:\cC\to \fdHilb$ is a fiber functor. 
    We recall that $\bbG$ is \emph{Kac-type} if and only if $F$ is a dimension-preserving functor. 
    
    From the cyclic module C*-category $\fdHilb$ with basepoint $\bbC,$ we obtain a connected C*-algebra object $\bfG$.
    By \cite[Proposition 6.4]{MR3948170}, $\bbG$ is Kac-type if and only if $\bfG$ is tracial. 
    Noteworthy, the class of Kac-type discrete quantum groups includes duals of all classical compact groups $\bbG = (\Rep(G), \mathsf{For})$, and all classical discrete groups. 
    Their corresponding connected C*-algebra objects $\bfG$ are therefore tracial. 

    For non-tracial examples, we consider the non-Kac compact quantum group $SU_q(2)$ for $q>0.$
    Letting $\delta = q+q^{-1}$, we have that $\Rep(SU_q(2))\cong \mathsf{TLJ}(\delta)$ is the \emph{Temperley-Lieb-Jones} unitary tensor category of parameter $\delta$ generated by a single strand. 
    The module categories for $\mathsf{TLJ}(\delta)$ correspond to \emph{fair and balanced} weighted graphs as explained in \cite{MR3420332}. 
    Considering pointed loops on these graphs, Jones and Penneys constructed connected C*-algebra objects, denoted $\bbA_v,$ in \cite[Definition 6.7]{MR3948170}, established necessary and sufficient conditions for traciality in \cite[Proposition 6.8]{MR3948170}, and described an explicit example from a weighted fair and balanced graph in \cite[Example 6.9]{MR3948170}.

    These conditions can be re-framed in terms of irreducible C*-discrete inclusions. 
    Given a discrete quantum group $\bbG = (\cC, F)$, there exists an outer action on a unital monotracial simple C*-algebra $A$ (c.f. \cite{MR4139893})
    $$H:\cC\to \fgpBim(A).$$
    We can then take reduced crossed product and obtain 
    $$
    A\overset{E}{\subset} A\rtimes_{H,r}\bfG,
    $$
    which is irreducible C*-discrete by design. 
    The corresponding connected C*-algebra object $\bfG$ is tracial if and only if $\bbG$ is Kac-type.     
\end{exs}

\subsection{Galois correspondence for outer actions of discrete groups}\label{sec:GalDisc}
We now recall Example \ref{ex:HilbgammOuterAction} and specialize our Galois correspondence to the case where $\cC=\fdHilb(\Gamma),$ where $\Gamma$ is a countable discrete group outerly acting on a simple unital C*-algebra $A$, and to re-express the Galois correspondence established in \cite{MR4010423} in categorified terms involving $\rCorr$.

Recall that for $x=\sum_{g\in\Gamma}x_g u_g\in A\rtimes_{\alpha,r}\Gamma$ we let the \textbf{support} of $x$ be $\supp(x)=\{g\in\Gamma|\ x_g=E(xu_g^{-1})\neq 0\},$ and for a subspace $X\subseteq A\rtimes_{\alpha,r}\Gamma,$ define $\supp(X):=\cup_{x\in X}\supp(x).$
Furthermore, we say $X$ is \textbf{replete} if it contains each $y\in A\rtimes_{\alpha,r}\Gamma$ for which $\supp(y)\subseteq\supp(X)$. This notion was studied by Cameron and Smith in \cite{MR4010423} under the name ``\emph{full},'' but we avoid this terminology because it conflicts with the already established notion of fullness for Hilbert C*-bimodules---where $M_A$ is full if the two-sided ideal of $A$ of inner products $\langle M\mid M\rangle_A$ is dense in $A$.

In \cite{MR4010423}, Cameron and Smith describe a bijection between the subsets $S\subseteq \Gamma$, and the $\|\cdot\|_B$-closed replete $A$-$A$ bimodules $A\subseteq Y\subseteq A\rtimes_{\alpha, r} \Gamma$ given by $S\mapsto Y_S:=\{y\in A\rtimes_{\alpha,r}\Gamma\mid \mathsf{supp}(y)\subseteq S\}.$ 
This bijection does not exactly fit into our categorical setting of right $A$-$A$ C*-correspondences, and the main goal of this section is to reformulate some of their results and re-express them using C* 2-categories leading to Proposition \ref{prop:FullvsProjective}.

\begin{lem}\label{lem:AbelianRelCom}
Let $A$ be a unital simple C*-algebra, $\Gamma$ be a countable discrete group, and $\alpha:\Gamma\to \Aut(A)$ be an outer action.
Then the relative commutant for the inclusion  $A\overset{E}{\subset}B:=A\rtimes_{r,\alpha}\Gamma = A\rtimes_{\alpha,r}\bbC[\Gamma]=:B$ is given by  
$$[A\rhd]'\cap \End^\dag(\cB_A) = \End^\dag({_A}\cB_A)\cong \ell^\infty(\Gamma).$$
Recall $E$ is the canonical conditional expectation. Furthermore, the inclusion is irreducible and C*-discrete.
\end{lem}
\begin{proof}
    Irreducibility for $A\subset B$ follows from the outerness of $\alpha$ by standard arguments.  
    Notice that the usual reduced crossed product of $A$ by $\Gamma$ by $\alpha$ coincides with the generalized reduced crossed product of $A$ by $\bbC[\Gamma]\in\rCorr_{\mathsf{fd}}(\fdHilb(\Gamma))$ by $\alpha$ as in \cite[\S3]{2023arXiv230505072H}. (See also \cite[Example 3.10]{2023arXiv230505072H}.)
    So $A\overset{E}{\subset} B\in\CDisc$ by \cite[Theorem 3.9]{2023arXiv230505072H}. 
    By Proposition \cite[Proposition 2.10]{2023arXiv230505072H} we know that $\End^\dag({_A}\cB_A)$ is $*$-isomorphic to an infinite sum of matrix algebras $M_{n_g}(\bbC),$ where each $n_g$ corresponds to the multiplicity of each isotypic irreducible bimodule indexed by $\Gamma=\Irr(\Vec(\Gamma))$. 
    But by the first part of the proof we know $n_g=1$ for each $g\in\Gamma,$ and so $\End^\dag({_A}\cB_A)$ is the abelian von Neumann algebra $\ell^\infty(\Gamma).$ 
\end{proof}
\begin{prop}\label{prop:FullvsProjective}
Let $A$ be a unital simple C*-algebra, $\Gamma$ be a countable discrete group, and $\alpha:\Gamma\to \Aut(A)$ be an outer action. 
Consider the associated C*-discrete inclusion $A\overset{E}{\subset} A\rtimes_{r,\alpha}\Gamma=:B.$ 
Then there is an equivalence 
\begin{align*}
    \left\{{}_AX_A\ \middle|\ 
        \begin{aligned}
            & X\subset B \text{ is a replete and}\\
            & \|\cdot\|_B\text{-closed } A\text{-}A \text{ bimodule }
        \end{aligned}
    \right\}
    \qquad\cong\qquad
    \left\{\cL\in\rCorr(A\to A)\ \middle|\ \cB \cong \cL\oplus \cL^\perp 
    \right\}.
    \end{align*}
\end{prop}
\begin{proof}
    Note that irreducibility and C*-discreteness for $A\overset{E}{\subset} B$ was established in Lemma \ref{lem:AbelianRelCom} above. 
    
    The correspondence is given as follows: For $X$ as above, consider the indicator function $\bbOne_{\supp(X)}\in\ell^\infty(\Gamma),$ which is a projection in $\End^\dag({_A}\cB_A)$ by Lemma \ref{lem:AbelianRelCom}. 
    We then map
    $$ X\mapsto\bbOne_{\supp(X)}[\cB]\in\rCorr(A\to A),$$
    and $\bbOne_{\supp(X)}[\cB]$ is an orthogonally complemented right $A$-$A$ C*-subcorrespondence of $\cB.$
    For the reversed mapping, let $\cL\cong \bbOne_L[\cB],$ where $L\subset \Gamma$ is the support of $\cL.$ 
    (This is well-defined by the Peter-Weyl decomposition from \cite[Corollary 2.9]{2023arXiv230505072H}.) 
    Using the notation from Example \ref{ex:HilbgammOuterAction}, let $\bbL(g):=\rCorr_{A-A}({_g}A\to\cL)^\diamondsuit:=\rCorr_{A-A}({_g}A\to\cL\cap B\Omega)= \rCorr_{A-A}({_g}A\to\cL) \cong \bbC\cdot\bbOne_{L}(g).$
    (So $\bbL\in\Vec(\fdHilb(\Gamma))$ directly corresponds to $\bbOne_L\in\ell^\infty(\Gamma)$.)
    Recall that C*-discreteness directly implies projective quasi-regularity, and so the equality above follows from \cite[Theorem 2.7]{2023arXiv230505072H}, saying every copy of ${_g}A$ contained in $\cL$ must lie inside $B\Omega.$ So, following the notation of \cite[\S3]{MR4010423} we map 
    $$\cL\mapsto Y_L:=\{b\in B|\ \supp(b)\subseteq L\}.$$
    (Recall, $\supp(b)=\{g\in\Gamma|\ E(g^{-1}\cdot b)\neq 0\}.$ )

    We shall now show these are mutual inverse mappings. 
    By \cite[Proposition 3.3 (ii)]{MR4010423}, since $X$ is assumed replete it follows that $X=Y_{\supp(X)},$ and so we get $X\mapsto \bbOne_{\supp(X)}[\cB]\mapsto Y_{\supp(X)}= X$ is the identity. 
    Furthermore, since $\cL$ is projective, we have $\supp(\bbL)=L$. 
    Thus, $\cL\mapsto Y_{L}\mapsto \bbOne_{\supp(Y_L)}[\cB] = \bbOne_L[\cB] = \cL$ also is the identity.  
\end{proof}

In the special case when $X=D$ is an intermediate C*-algebra, we can now rephrase \cite[Theorem 3.5]{MR4010423} (corresponding the intermediate subalgebras $A\subset D\subset B$ with subgroups $\Lambda\leq\Gamma$) in the language of our Corollary \ref{cor:Galois}. 
We remark that Cameron and Smith showed that any such $D$ is the target of a surjective conditional expectation $E^B_D:B\to D$, automatically making the inclusion $A\overset{E|_D}{\subset}D$ C*-discrete and implying that ${_A}D_A$ is replete, 
whereas in Theorem \ref{thm:ceCDisc} we characterized intermediate C*-discrete inclusions precisely as the targets of compatible conditional expectations. 
So in this case, the lattice of intermediate C*-algebras is exhausted by the connected C*-algebra objects.
\begin{cor}
    Let $A$ be a unital simple C*-algebra, $\Gamma$ be a countable discrete group, and $\alpha:\Gamma\to \Aut(A)$ be an outer action, 
    and consider the inclusion  $A\overset{E}{\subset}A\rtimes_{\alpha,r}\Gamma = A\rtimes_{\alpha,r}\bbC[\Gamma]=:B.$ 
    Then there is a correspondence    
    \begin{align*}
    &\left\{ D \ \middle|\ 
        \begin{aligned}
            & A\subseteq D\subseteq B \text{ is unital }\\
            & \text{ intermediate  C*-algebra }
        \end{aligned}
    \right\}
    \ \ \cong&&
    \left\{\bbD\in\rCorr(\fdHilb(\Gamma))\ \middle|\ \text{ and } \bbOne\subseteq\bbD\subseteq \bbB  \right\}.
    \end{align*}
\end{cor}
\begin{proof}
    As we have seen above, the objects in $\Vec(\fdHilb(\Gamma))$ can be interchanged for projections in $\ell^\infty(\Gamma),$ which in this case necessarily contain the group identity $e\in\Gamma$ inside its support. 
    Furthermore, the subgroup condition on the support clearly corresponds to algebra objects.
\end{proof}
\begin{remark}
    We thank Prof. Roger Smith for further clarifying \cite[Theorem 3.5]{MR4010423} implying that intermediate C*-algebras $A\subseteq D\subseteq A\rtimes_{\alpha,r}\Gamma =:B$ are necessarily ``full''/replete. This is, setting $\mathsf{supp}(D)=:\Lambda\leq \Gamma,$ gives  $\overline{\mathsf{span}\{x_hu_h|h\in\Lambda, x_h\in A\}}^{\|\cdot\|_B} = D = \{y\in B|\ \mathsf{supp}(y)\subseteq \Lambda\}$ because intermediate C*-algebras are targets of compatible conditional expectations. 
    Thus, our setting yields a correct categorification of Cameron-Smith's bimodule formalism.
\end{remark}

\section{Actions of UTCs on C*-algebras}

\subsection{Free actions of UTCs}
In this section we propose a definition of \emph{freeness} for an action $F:\cC\to\fgpBim(A)$, where $\cC$ is a UTC and $A$ is a unital C*-algebra. We will then crucially use this property to ensure that if $A$ is simple, $F$ is outer, and $\bbB\in\Vec(\cC^{\op})$ is a connected $\cC$-graded C*-algebra, then the reduced crossed product C*-algebra of $A$ by $F$ under $\bbB$ denoted $B:=A\rtimes_{F,r}\bbB$ remains simple. 

\begin{defn}\label{defn:FreeAction}
    Let $\cC$ be a unitary tensor category, and let $A$ be a unital C*-algebra. 
    We say that an outer action (i.e. fully-faithful unitary tensor functor) $F:\cC\to \fgpBim(A)$ is \textbf{free} if for each $c\in\cC$ not containing $1_\cC$ in its semisimple decomposition and every $\xi\in F(c)$ we have 
    $$0=\inf\left\{\left\|\sum_{i=1}^n a_i^*\rhd\xi\lhd a_i\right\|\ \Bigg|\ a_1,..., a_n\in A,\ \sum_{i=1}^n a_i^*a_i =1\right\}.$$
\end{defn}

This condition stems from categorifying the averaging properties for outer actions of discrete groups on unital C*-algebras obtained by \cite{MR4010423}. There, they show that for discrete group actions freeness is a consequence of outerness. One should also compare this freeness condition for UTC-actions with \emph{Popa's averaging technique} for collections of outer automorphisms of C*-algebras with the \emph{Dixmier Property} \cite[Lemma 5.10]{MR4599249}. Definition \ref{defn:FreeAction} should also be compared with \emph{Kishimoto's condition} \cite[Definition 2.8]{2016arXiv161106954K}, employed to \emph{detect} ideals in $B$ with ideals in $A$, and first appearing in \cite{MR634163}.

We now show that connected C*-discrete extensions of a unital simple C*-algebra by a free quantum dynamics are also simple by directly adapting the techniques from \cite[\S3]{MR4010423}.  
\begin{thm}[{Theorem~\ref{thmalpha:CrossProdSimple}}]\label{thm:CrossProdSimple}
    Let $\cC$ be a UTC, and $A$ be a unital simple C*-algebra with a free and outer action $F:\cC\to \fgpBim(A).$ 
    Let $\bbB\in\Vec(\cC^{\op})$ be a connected $\cC$-graded C*-algebra and consider the generalized reduced crossed product inclusion $A\overset{E}{\subset} A\rtimes_{r,F}\bbB:= B,$ where $E$ is the canonical conditional expectation onto $A$. 
    Then $B$ is unital and simple.
\end{thm}
\begin{proof}
    Let $I\subseteq B$ be a nonzero closed two-sided ideal. 
    For $c\in\cC$ we let $P_c\in\End^\dag({_A}\cB_A)$ be the orthogonal projection onto the isotypic component of $F(c)$ in $\cB,$ which is explicitly given by a Pimsner--Popa basis. 
    For $b\in B,$ let $$\supp_\cC(b):=\{c\in\Irr(\cC)\ |\ P_c(b\Omega)\neq 0\},$$ and for $X\subset B$ we let $\supp_\cC(X):=\cup_{x\in X}\supp_\cC(x).$ 
    Let $x_0\in I$ be nonzero, and by replacing with $x_0^*x_0\in I$ without loss of generality, we may assume that $1_\cC\in\supp_\cC(x_0).$ 
    (Since $E$ is faithful, this coefficient is $E(x_0^*x_0)\neq 0$.) 
    Since $A$ is simple and unital, it follows that $1\in AE(x_0)A =\{\sum_{i\in I} a_iE(x_0)a_i'\ |\ I \text{ is finite} \},$ so without loss of generality we can further assume that $E(x_0)=1$. 
  
    Denote $\PQN(A\subset B)$ by $B^{\diamondsuit}$, and let $\epsilon\in (0,1)$ and 
    $$(I\times B^{\diamondsuit})_\epsilon:=\{(x,y)\ |\ x\in I, E(x)=1, y\in B^{\diamondsuit}, \|x-y\|<\epsilon\}.$$
    Since $A\subset B$ is C*-discrete by \cite[Theorem 3.9]{2023arXiv230505072H}, then $\overline{B^{\diamondsuit}}^{\|\cdot\|_B} = B$, and so for some $y\in B^{\diamondsuit}$ we obtain $(x_0,y)\in (I\times B^{\diamondsuit})_\epsilon \neq\emptyset.$ 
    
    Let $(x,y)\in (I\times B^{\diamondsuit})_\epsilon$ be fixed but arbitrary with the property that $y=\sum_{c\in \Irr(\cC)}\xi_{(c)}\otimes y_{(c)}$ is an element of minimal length considered with $\Irr(\cC)$ as endowed with the partial order given by containment.  We then have that $\|1-E(y)\|=\|E(x-y)\|\leq \|x-y\|<\epsilon <1.$ 
    Thus $E(y)\neq 0.$ By means of contradiction, we shall now show that $\supp_\cC(y)$ is the singleton $\{1_\cC\}$ and hence $y\in A$. Assume there is some $c\in\supp_\cC(y)\setminus\{1_\cC\}{}$. Let $\delta>0$ be such that $\|x-y\|+\delta<\epsilon.$ Now, since $F$ was assumed free, let $\{a_i\}_1^n\subset A$ with $\sum_{i=1}^n a_i^* a_i=1$ and such that the bound $\|\sum_{i=1}^n a_i^*\rhd P_c(y\Omega)\lhd a_i\|<\delta$ holds.
    
    Now consider the contraction 
    \begin{align*}
        \Psi: B&\to B\\
        b&\mapsto \sum_1^n a_i^*ba_i,
    \end{align*}
    and we notice at once that: $\Psi[I]\subset I,$ $E(\Psi(x))=1$, and $\|P_c(\Psi(y)\Omega)\|<\delta.$ Also $A B^{\diamondsuit} A\subset B^{\diamondsuit}$ and $P_c[\cB]\subset B^{\diamondsuit} \Omega$ implies $\Psi(y)-\check{P_c}(\Psi(y)\Omega)\in B^{\diamondsuit}$. 
    Using the Pimsner--Popa inequality from \cite[Lemma 3.6]{2023arXiv230505072H} and shrinking $\delta$ if necessary, we obtain $\| \check{P_c}(\Psi(y)\Omega)\| <\delta$. Thus
    $$\|\Psi(x)-(\Psi(y)-\check{P_c}(\Psi(y)\Omega))\| \leq \|\Psi(x-y)\|+\delta \leq \|x-y\|+\delta < \epsilon.$$ 
    Hence $(\Psi(x), \Psi(y)-\check{P_c}(\Psi(y)\Omega))\in (I\times B^{\diamondsuit})_\epsilon$. But $\supp_\cC\left(\Psi(y)-\check{P_c}(\Psi(y)\Omega)\right)$ is strictly smaller than $\supp_\cC(y);$ a contradiction. 
    Therefore, the assumption that $\supp_\cC(y)\neq \{1_\cC\}$ is absurd, and so $y\in A.$
     
     Repeating the same argument with $\epsilon=1/2^n$ for $n\geq 1$ yields sequences $\{x_n\}_n\subset I$ and $\{y_n\}_n\subset A$ with $(x_n,y_n)\in (I\times B^{\diamondsuit})_{2^{-n}}.$ 
     But since $\|1-y_n\| = \|E(x_n-y_n)\|<2^{-n}$ it then follows that $\|1-y_n\|\to 0 $ in $\|\cdot\|_B$-norm and therefore $\|x_n-1\|\to 0$ as well. 
     Then $1\in I$ and finally $I=B$, so $B$ is simple. 
\end{proof}

\section{\texorpdfstring{$A$}{A}-Valued Semicircular Systems}\label{sec:A-valued_semicirc}
Let $A$ be a unital C*-algebra admitting a family of linear maps $\eta_{ij}\colon A \to A$ for $i,j$ in a countable set $I$. One can package these maps into a single map
    \begin{align*}
        \eta\colon A&\to A\otimes B(\ell^2(I))\\
                    a&\mapsto \sum_{i,j\in I} \eta_{ij}(a)\otimes e_{ij},
    \end{align*}
where $e_{ij}\in B(\ell^2(I))$ are matrix units. One says $\eta$ is a \emph{covariance matrix} if it is completely positive. In \cite{Shl99}, Shlyakhtenko considers the case when $A$ is a von Neumann algebra and each $\eta_{ij}$ is a normal completely positive map, and constructs a von Neumann algebra $\Phi(A,\eta)$ generated by $A$ and a so-called \emph{$A$-valued semircular family} $\{X_i\colon i\in I\}$ \emph{with covariance $\eta$}. This construction is astonishingly versatile in that it can yield virtually every von Neumann algebra of interest in free probability theory: (interpolated) free group factors, free Araki--Woods factors, amalgamated free products, etc. This construction can be restricted to the category of C*-algebras (see \cite[Section 2.10]{Shl99}), and we will show in this section that it yields C*-discrete inclusions under mild conditions on $\eta_{ii}$, $i\in I$ (see Theorem~\ref{thm:A-valued_discrete_inclusions}).

Before proceeding, recall that for a completely positive map $\theta\colon A\to B$ between C*-algebras one can define a right $A$-$B$ C*-correspondence through separation and completion of $A\odot B$ with respect to the right $B$-valued inner product:
    \[
        \<a_1\otimes b_1\mid a_2 \otimes b_2\>_B:= b_1^* \theta(a_1^* a_2)b_2 \qquad a_1,a_2\in A,\ b_1,b_2\in B.
    \]
(See also \cite[Chapter 5]{MR1325694}). We  denote this C*-correspondence by $A\otimes_\theta B$. Note that if $A$ is unital and $\lambda$ is a unital $*$-endomorphism, then one has $A\otimes_\lambda A \cong {}_\lambda A_A$ and $A\otimes_\lambda A$ is \emph{algebraically} generated by $[1\otimes 1]$. Additionally,
    \begin{align}\label{eqn:tensor_composition}
        (A\otimes_\lambda A)\boxtimes_A (A\otimes_\rho A) \cong A\otimes_{\rho\circ \lambda} A \cong {}_{\rho\circ\lambda} A_A,
    \end{align}
for $\rho$ a unital $*$-endomorphism.

Let $A$ be a unital C*-algebra and  let $\eta$ be a covariance matrix. Fix $i_0\in I$ and define
    \[
        T:=A\otimes_\eta (A\otimes B(\ell^2(I)))\lhd(1\otimes e_{i_0i_0}),
    \]
which is a right $A$-$A$ C*-correspondence with the actions it inherits from $A\otimes_\eta (A\otimes B(\ell^2(I)))$ and the right $A$-valued inner product
    \[
        \<\xi\mid \eta\>_A:= (1\otimes \Tr)\left(\<\xi\mid \eta\>_{A\otimes B(\ell^2(I))}\right).
    \]
Moreover, $T$ is topologically generated by $\{\xi_i\colon i\in I\}$ where
    \[
        \xi_i := [1\otimes (1\otimes e_{i i_0})].
    \]
Observe that for $a\in A$ and $i,j\in I$ one has
    \begin{align}\label{eqn:eta_ij_inner_product}
    \begin{aligned}
        \< \xi_i\mid  a\rhd\xi_j\>_A &=\< [1\otimes (1\otimes e_{i i_0})], [a\otimes (1\otimes e_{j i_0})]\>_A \\    &=  (1\otimes \Tr)\left((1\otimes e_{i_0i})\eta(a)(1\otimes e_{j i_0})\right) = \eta_{ij}(a).
    \end{aligned}
    \end{align}
Also note that $T$ does not depend on the choice of $i_0\in I$. Consider the Fock space associated to $T$:
    \[
        \cF(T):= \overline{\bigoplus_{d\geq 0}}\  T^{\boxtimes_A d},
    \]
where $T^{\boxtimes_A 0}:=A$. Let $\hat{1}\in T^{\boxtimes_A 0}$ denote the copy of $1\in A$, so that $A\hat{1} = T^{\boxtimes_A 0}$. For each $\xi\in T$, define $\ell(\xi)\in \End^\dag(\cF(T)_A)$ in the usual way:
    \begin{align*}
        \ell(\xi) a\hat{1} &= \xi \lhd a\\
        \ell(\xi) \eta_1\otimes \cdots \otimes \eta_d &= \xi\otimes \eta_1 \otimes \cdots \otimes \eta_d.
    \end{align*}
The adjoint is given by
    \begin{align*}
        \ell(\xi)^* a \hat{1} &= 0\\
        \ell(\xi)^* \eta_1\otimes \cdots \otimes \eta_d &= \<\xi\mid \eta_1\>_A \rhd \eta_2 \otimes \cdots \otimes \eta_d.
    \end{align*}
The (diagonal) left action of $A$ allows us to identify $A\subset \End^\dag(\cF(T)_A)$, and one has $a\ell(\xi)b = \ell(a\rhd \xi \lhd b)$ for all $a,b\in A$ and $\xi\in T$. For each $i\in I$, define
    \[
        X_i:=\ell(\xi_i) + \ell(\xi_i)^*,
    \]
and define
    \[
        \hat{\Phi}(A,\eta):=C^*(A\cup \{X_i\colon i\in I\}) \subset \End^\dag(\cF(T)_A).
    \]
The map $E\colon \hat{\Phi}(A,\eta)\to A$ defined by
    \[
        E(x):=\<\hat{1}\mid x\hat{1}\>_A
    \]
defines a conditional expectation, which by (\ref{eqn:eta_ij_inner_product}) satisfies
    \[
        E(X_i a X_j) = \< \xi_i \mid a \rhd \xi_j\>_A = \eta_{ij}(a)
    \]
for all $i,j\in I$ and $a\in A$. In general, $E$ need not be faithful, but if $A$ admits a faithful state (e.g. if $A$ is separable) then \cite[Proposition 5.2]{Shl99} characterizes when $E$ is faithful. In this case the right $\hat{\Phi}(A,\eta)$-$A$ C*-correspondence associated to $E$ has a very explicit representation:

\begin{prop}\label{prop:A-valued_standard_form}
Let $A$ be a unital C*-algebra with covariance matrix $\eta\colon A\to A\otimes B(\ell^2(I))$ such that the conditional expectation $E\colon \hat{\Phi}(A,\eta)\to A$ is faithful. Denote by $\cB:=\overline{\hat{\Phi}(A,\eta)\Omega}^{\|\cdot\|_A}$ the right $\hat{\Phi}(A,\eta)$-$A$ C*-correspondence with $\<x\Omega\mid y\Omega\>_A:=E(x^*y)$. Then $x\Omega\mapsto x\hat{1}$ extends to an isometric isomorphism
    \[
        {}_{\hat{\Phi}(A,\eta)} \cB_{A}\cong {}_{\hat{\Phi}(A,\eta)} \cF(T)_{A}.
    \]
In particular, for all $i\in I$ one has
    \[
        \overline{A\rhd X_i\Omega \lhd A} \cong \overline{A\rhd \xi_i \lhd A} \cong A\otimes_{\eta_{ii}} A
    \]
as right $A$-$A$ C*-correspondences.
\end{prop}
\begin{proof}
By definition of $E$ we have
    \[
        \<x\Omega\mid y\Omega\>_A = E(x^*y) = \<x\hat{1} \mid y\hat{1}\>_A,
    \]
for all $x,y\in \hat{\Phi}(A,\eta)$. So it suffices to show $\hat{\Phi}(A,\eta)\hat{1}$ is dense in in $\cF(T)$.  For each $d\geq 1$ denote
    \[
        S^{(d)}:= \text{span}\{\eta_1\boxtimes \cdots \boxtimes \eta_d \in T^{\boxtimes_A d}\colon \eta_j\in A\rhd \xi_{i_j} \lhd A,\ i_1,\ldots i_d\in I\},
    \]
and $S^{(0)}:=A\hat{1}$. For each $d\geq 0$ and each $\eta\in S^{(d)}$, we claim there exists $x\in \bbC\<A\cup \{X_i\colon i\in I\}\>$ satisfying $x\hat{1}=\eta$. This is immediate for $d=0$, so we suppose it has been shown for $0,1,\ldots, d-1$. Consider
    \begin{align}\label{eqn:element_of_spanning_set_of_Sd}
      \zeta:=(a_1\rhd \xi_{i_1} \lhd b_1)\boxtimes \cdots \boxtimes (a_d\rhd \xi_{i_d} \lhd b_d) \in S^{(d)}  
    \end{align}
for some $a_1,b_1,\ldots, a_d,b_d\in A$ and $i_1,\ldots, i_d\in I$. Then by definition of the $X_i$ we have
    \[
        (a_1 X_{i_1} b_1 \cdots a_d X_{i_d} b_d)\hat{1} - \zeta \in \text{span}\left( S^{(0)}\cup \cdots \cup S^{(d-1)}\right).
    \]
Our induction hypothesis allows us to find $x_0 \in \bbC\<A\cup \{X_i\colon i\in I\}\>$ so that $x_0\hat{1}$ is the above difference. Taking $x:= a_1 X_{i_1} b_1 \cdots a_d X_{i_d} b_d -x_0$ gives $x\hat{1} = \zeta$. Since elements of the form in (\ref{eqn:element_of_spanning_set_of_Sd}) span $S^{(d)}$, we have established the claim. Thus
    \[
        \hat{\Phi}(A,\eta) \hat{1} \supset \bigoplus_{d\geq 0} S^{(d)}
    \]
is dense in $\cF(T)$.

The first computation in this proof gives the first isomorphism in the final claim. To see the second isomorphism, we use (\ref{eqn:eta_ij_inner_product}) to compute
    \[
        \< a \rhd \xi_i \lhd b \mid c \rhd \xi_i \lhd d \>_A = b^* \eta_{ii}(a^*c) d = \< [a\otimes b]\mid [c\otimes d]\>_{A},
    \]
for all $a,b,c,d\in A$ and $i\in I$.
\end{proof}

We now give some sufficient conditions for the inclusion $A\overset{E}{\subset} \hat{\Phi}(A,\eta)$ to be C*-discrete in the case that $A$ admits a faithful tracial state.

\begin{thm}\label{thm:A-valued_discrete_inclusions} 
Let $A$ be a unital C*-algebra with covariance matrix $\eta\colon A\to A\otimes B(\ell^2(I))$. Suppose that:
    \begin{enumerate}[label=(\roman*)]
    \item $A$ admits a faithful tracial state $\tau$ satisfying $\tau(\eta_{ij}(x) y) = \tau(x\eta_{ji}(y))$ for all $i,j\in I$ and $x,y\in A$; and

    \item $A\otimes_{\eta_{ii}} A\in \fgpBim(A)$ for all $i\in I$.
    
    \end{enumerate}
Then $\tau\circ E$ is a faithful tracial state on $\hat{\Phi}(A,\eta)$ and $A\overset{E}{\subset} \hat{\Phi}(A,\eta)$ is a C*-discrete inclusion. In the case that $I$ is finite, the inclusion is irreducible if and only if the $A$-central vectors in the Fock space are $\cF(T)^A = \bbC\hat{1}$.
\end{thm}
\begin{proof}
Condition (i) and \cite[Proposition 2.20]{Shl99} imply $\tau\circ E$ is tracial. Moreover, for $a,b,c,d\in A$ and $i,j\in I$ one has
    \begin{align*}
        \tau(\<a\rhd \xi_i \lhd b\mid c\rhd \xi_j \lhd d \>_A) &= \tau( b^*\eta_{ij}(a^*c)d)\\
        &= \tau( c \eta_{ji}(db^*)a^*) = \tau(\<d^*\rhd \xi_j \lhd c^* \mid b^* \rhd \xi_i \lhd a^*\>_A).
    \end{align*}
This shows the map $J\colon a\rhd \xi_i \lhd b\mapsto b^* \rhd \xi_i \lhd a^*$ extends to a conjugate linear isometry on the completion of $T$ with respect to the $\bbC$-valued inner product $\tau(\<\ \cdot\mid -\>_A)$. Thus $E$ (and hence $\tau\circ E$) is faithful by \cite[Propsition 5.2]{Shl99}.  

Now, from Proposition~\ref{prop:A-valued_standard_form} we have for all $i\in I$
    \[
        \overline{A\rhd X_i\Omega \lhd A} \cong \overline{A\rhd \xi_i \lhd A} \cong A\otimes_{\eta_{ii}} A,
    \]
and by condition (ii) and Theorem ~\ref{thm:ProjectiveVectors} the latter is algebraically generated by $[1\otimes 1]$. This implies $$\overline{A\rhd X_i\Omega \lhd A} = A\rhd X_i\Omega \lhd A \subset \hat{\Phi}(A,\eta)\Omega.$$ Thus $X_i\in \PQN(A\subset \hat{\Phi}(A,\eta))$ for all $i\in I$.

Finally, suppose $I$ is finite so that we can assume $I=\{1,\ldots, n\}$ after relabeling if necessary. We will show
    \[
        \cF(T)^A = \overline{[A'\cap \hat{\Phi}(A,\eta)]\hat{1}},
    \]
which gives the final claim since $\hat{1}$ is separating for $\hat{\Phi}(A,\eta)$ by virtue of $E$ being faithful. The inclusion `$\supset$' is clear, so suppose $\gamma\in \cF(T)$ is $A$-central. By definition of the Fock space, we can orthogonally decompose $\gamma = \sum_{d\geq 0} \gamma_d$ where $\gamma_d\in T^{\boxtimes_A d}$ is also $A$-central. Now, condition (ii) and \cite[Lemma 1.8 (2)]{2023arXiv230505072H} imply $T\in \fgpBim(A)$ since it is the image of the bounded $A$-$A$ bimodular map
    \[
        \bigoplus_{i=1}^n (A\rhd \xi_i \lhd A) \ni \zeta_1\oplus \cdots \oplus \zeta_n\mapsto \zeta_1+\cdots +\zeta_n\in T.
    \]
Consequently, $T^{\boxtimes_A d} \in \fgpBim(A)$ for all $d\geq 0$. In particular, each  $T^{\boxtimes_A d}$ equals $S^{(d)}$ from the proof of Proposition~\ref{prop:A-valued_standard_form}. Indeed, $T^{\boxtimes_A d}$ is algebraically generated by some right Pimsner--Popa basis, which can be approximated by elements in $S^{(d)}$, and therefore $S^{(d)}$ contains algebraic generators by Mostow's Lemma~\ref{lem:BasisApproximation}. Checking the proof of Proposition~\ref{prop:A-valued_standard_form} again, we see 
 that each $\gamma_d$ is of the form $\gamma_d=y_d\hat{1}$ for some $y_d\in \hat{\Phi}(A,\eta)$. For all $x\in A$ we then have
    \[
        \|E( |xy_d-y_dx|^2)\| = \|(xy_d-y_dx)\hat{1}\|_A^2 = \| x\rhd \gamma_d - \gamma_d\lhd x \|_A^2 = 0.
    \]
That is, $y_d\in A'\cap \hat{\Phi}(A,\eta)$ for all $d\geq 0$ and therefore $\gamma$ is in the claimed closure.
\end{proof}

\bibliographystyle{amsalpha}
\bibliography{bibliography}
\end{document}